\documentclass[preprint,12pt,authoryear]{elsarticle}
\usepackage{setspace}

\usepackage{comment}
\usepackage{amssymb}
\usepackage{amsmath}
\allowdisplaybreaks 
\usepackage[ruled,vlined,linesnumbered]{algorithm2e}
\usepackage{booktabs}
\usepackage{threeparttable}
\usepackage{enumitem}
\usepackage{pgfplots}
\usepackage{pgfplotstable}
\usepackage{filecontents}
\pgfplotsset{compat=1.17}
\usepackage{tikz}
\usepackage{xurl}
\usepackage{longtable}
\usepackage{array}
\usepackage{geometry}

\usepackage[symbol]{footmisc}

\journal{}

\begin{document}
\onehalfspacing
\begin{frontmatter}



\title{Scheduling Electricity Production Units to Mitigate Severe Weather Impact: An Efficient Computational Implementation} 

 \author[label1]{Yongzheng Dai}
 \author[label2,label3]{Antonio J. Conejo}
 \author[label4]{Feng Qiu}
  \affiliation[label1]{organization={Georgia Institute of Technology, ISyE},
 	addressline={755 Ferst Drive N.W.},
 	city={Atlanta},
 	postcode={30332},
 	state={GA},
 	country={USA}}
 \affiliation[label2]{organization={The Ohio State University, ISE},
             addressline={1971 Neil Avenue},
             city={Columbus},
             postcode={43212},
             state={OH},
             country={USA}}
 \affiliation[label3]{organization={The Ohio State University, ECE},
	addressline={2015 Neil Avenue},
	city={Columbus},
	postcode={43212},
	state={OH},
	country={USA}}
 \affiliation[label4]{organization={Argonne National Laboratory},
             addressline={9700 S. Cass Avenue},
             city={Lemont},
             postcode={60439},
             state={IL},
             country={USA}}



\begin{abstract}
{
Extreme weather events in electric power systems can cause line trips or physical damage to transmission infrastructure, potentially leading to large-scale load shedding. To mitigate this risk, we propose a framework that strategically pre-positions the commitment of generation units—particularly slow-start units—to adapt to transmission topologies that may arise following such events. The objective is to minimize load shedding under worst-case conditions.
This paper makes two main contributions. First, we provide a more accurate representation of the underlying physical laws than those used in prior studies. Second, we develop a highly efficient solution algorithm that outperforms state-of-the-art, off-the-shelf solvers.
The proposed framework is formulated as a two-stage robust optimization model. In the first stage, generation units are scheduled in anticipation of disruptions. In the second stage, power dispatch decisions are optimized to minimize load shedding under the worst-case transmission topology. To ensure system reliability and security, we incorporate convexified AC power flow constraints. The resulting model is a tri-level mixed-integer nonlinear optimization problem.
To address its computational complexity, we design a problem-specific outer approximation algorithm embedded within a column-and-constraint generation framework. Computational results show that the proposed model and solution approach can achieve solutions within a standard optimality gap in reasonable time for moderately large instances.
}
\end{abstract}





\end{frontmatter}



\section*{Nomenclature}

\begin{longtable}{>{\bfseries}p{0.2\textwidth} p{0.7\textwidth}}
	\multicolumn{2}{l}{(1) Sets and Indexes:}\\
	$\mathcal{T}$ &Time periods, $\mathcal{T}:=\{1,...,T\}$.\\
	$\mathcal{G}$ &Generating units, $\mathcal{G}:=\mathcal{G}^\mathrm{T}\cup \mathcal{G}^\mathrm{R}$.\\
	$\mathcal{G}^\mathrm{T}$ &Thermal generating units.\\
	$\mathcal{G}^\mathrm{R}$ &Weather-dependent renewable generating units.\\
	$\mathcal{N}$ &Nodes in the electric network.\\
	$\mathcal{R}$ &Reliability Areas.\\
	$\Omega_{r}$ &Generating units in the reliability area $r$.\\
	$\mathcal{E}$ &Transmission lines.\\
	$\Lambda_n$ &Adjacent nodes to node $n$.\\
	$\Lambda_n^\mathrm{G}$ &Generating units at node $n$.\\
	$\Omega_{\mathrm{traj}}$ &Line status under all hurricane trajectories.\\
	\multicolumn{2}{l}{(2) Constants:}\\
	$C^\mathrm{F}_g$ &Fixed cost of unit $g$.\\
	$C^\mathrm{SU}_g/C_g^\mathrm{SD}$ &Startup/shutdown cost of unit $g$.\\
	$C_g^\mathrm{V}$ &Variable cost of unit $g$.\\
	$C^\mathrm{U}$ &Unserved energy cost.\\
	$L_g/F_g$ &Initial up/down time of thermal generating unit $g$.\\
	$T_g^\mathrm{U}/T_g^\mathrm{D}$ &Minimum up/down time of thermal generating unit $g$.\\
	$p_g^{\min}/p_g^{\max}$ &Minimum/maximum active power output from unit $g$.\\
	$q_g^{\min}/q_g^{\max}$ &Minimum/maximum reactive power output from unit $g$.\\
	$u_{g,0},y_{g,0},z_{g,0}$ &Status at the beginning of the scheduling horizon of unit $g$.\\
	$p_{g,0}$ &Initial active power output of unit $g$.\\
	$R_g^\mathrm{U}/R_g^\mathrm{SD}$ &Up/startup ramping limit of unit $g$.\\
	$R_g^\mathrm{D}/R_g^\mathrm{SD}$ &Down/shutdown ramping limit of unit $g$.\\
	$R_{r,t}^\mathrm{D}$ &Reserve required in reliability area $r$ in period $t$.\\
	$S_{n,m}$ &Capacity of line $(m,n)$.\\
	$G_{n,m}/B_{n,m}$ &Conductance/susceptance of line $(m,n)$.\\
	$b_{n,m}^\mathrm{shunt}$ &Half of the shunt susceptance of line $(m, n)$.\\
	$V_n^{\max}/ V_n^{\min}$ &Maximum/minimum voltage level for node $n$.\\
	$p^\mathrm{D}_{n,t}/ q_{n,t}^\mathrm{D}$ &Active/reactive power demand in period $t$ and node $n$.\\
	$n_\mathrm{traj}$ &Number of alternative hurricane trajectories.\\
	$a_{n,m}^{(k)}$ &On/off status of line $(n,m)$ under the $k$-th hurricane trajectory.\\
	$T$ &Number of stages.\\
	\multicolumn{2}{l}{(3) Variables:}\\
	$u_{g,t}$ &On/off status in period $t$ of unit $g$.\\
	$y_{g,t}$ &Startup indicator in period $t$ of unit $g$.\\
	$z_{g,t}$ &Shutdown indicator in period $t$ of unit $g$.\\
	$p_{g,t}/q_{g,t}$ &Active/reactive power produced by unit $g$ in period $t$.\\
	$\bar{p}_{g,t}$ &Maximum available active power output in period $t$ from unit $g$.\\
	$p_{n,m,t}/q_{n,m,t}$ &Active/reactive power flow of line $(n,m)$ in period $t$.\\
	$p_{n,t}^\mathrm{U}/q_{n,t}^\mathrm{U}$ &Unserved active/reactive load in period $t$ and node $n$.\\
	$f_{n,m,t}$ & Auxiliary variable equal to $S_{n,m,t}$.\\
	$a_{n,m}$ &On/off status of line $(n,m)$; $a_{n,m}=1$ if on and $a_{n,m} = 0$ if off.\\
	$c_{n,m,t}, s_{n,m,t}$ &Auxiliary variables for the second-order-conic formulation.\\
	$c_{n,n,t}^{m}$ &Auxiliary variable to linearize $a_{n,m}c_{n,n,t}$.\\
    $D_{n,m,t}^{(k)}$ &4 auxiliary variables ($k=1,2,3,4$) of second-order-conic expression for line $(n,m)$ in period $t$.\\
	$\omega_k$ &Hurricane trajectory indicator.
\end{longtable}

\section{Introduction}

{
	\subsection{Motivation, formulation, and solution technique}

Extreme weather events result in substantial economic losses worldwide. \cite{campbell2012weather} indicates that major power outages, primarily caused by storm-related damage to electricity transmission lines, cost an estimated \${20}--\${55} billion annually due to lost orders, spoiled inventory, delayed production, and other disruptions. On a global scale, extreme weather events—including hurricanes, floods, droughts, and storms—have caused more than \${2} trillion in damages over the past decade, with the United States accounting for nearly half of this total~\citep{ICCOxera2024}.

Severe weather also inflicts significant damage on energy infrastructure. For example, hurricanes Katrina and Rita damaged or destroyed thousands of miles of transmission lines and hundreds of substations; in Southeast Louisiana and Mississippi, 1.1 million customers lost power after Katrina~\citep{us_senate2005katrina}. In 2017, Hurricane Maria caused power outages affecting nearly all customers in Puerto Rico as of September 20, 2017~\citep{DOE2017hurricanes}. More recently, Hurricane Ida (2021) led to the failure of all eight transmission lines supplying New Orleans, cutting power to nearly 895{,}000 customers in Southeastern Louisiana~\citep{ENR2021ida}. These examples illustrate how transmission line failures due to wind, flooding, or structural collapse can cascade throughout the system, requiring emergency redispatch of generation or involuntary load shedding to maintain system stability.

To mitigate the adverse impacts of extreme weather events, various strategies have been studied and implemented, including infrastructure hardening, operational improvements \citep{abdelmalak2022proactive,DOE_PR_Resilience_2018}, and robust planning \citep{sang2019integrated,mohammadi2020tractable}. Significant investments have been made to strengthen infrastructure; for example, Florida Power and Light (FPL) has invested approximately \$3 billion in grid hardening efforts~\citep{FPL2017restoration}. At the operational level, researchers have explored both proactive and corrective control actions, such as adjusting generator commitments and dispatch decisions in anticipation of component outages \citep{wang2016resilience,abdelmalak2022proactive}.

In this paper, we develop an AC-based robust generation scheduling framework to mitigate the impact of severe weather events on transmission systems. Our primary objective is to ensure robust system security under extreme conditions, which leads to an adaptive robust optimization formulation with three interrelated levels. This approach is particularly important given the increasing frequency and severity of such events. We model the scheduling of generation units (the unit commitment problem) while representing the transmission network using AC power flow constraints in second-order conic form. As is customary, we apply a slight convex relaxation to this formulation, which may introduce minor inaccuracies in flow and generation levels but typically does not affect commitment decisions.

The first level of the three-level adaptive robust optimization model determines the unit commitment decisions that minimize cost while anticipating the worst-case weather scenario. Given these decisions, the second-level problem identifies the most disruptive weather event, maximizing system cost under the assumption that the system operator responds optimally. Finally, the third-level problem represents system operation under the realized worst-case event, minimizing operational cost given the first-stage commitment decisions. Overall, this results in a challenging min–max–min optimization problem with second-order conic constraints.

Since the third-level problem corresponds to a multi-period optimal power flow model that is convex and second-order conic, it can be replaced exactly by its dual formulation. This allows it to be merged with the second-level problem, yielding a two-level reformulation that can be solved using the column-and-constraint generation algorithm introduced by \cite{Bo2017Solving}.

The column-and-constraint generation approach involves a master problem and a subproblem. The master problem is a mixed-integer second-order conic program whose size grows significantly with the number of iterations, making it the primary computational bottleneck. The subproblem, in contrast, is a continuous second-order conic program of moderate size. To alleviate the computational burden of the master problem, we employ an effective problem-specific outer approximation technique.
}
\subsection{Literature Review}

Second-order-conic (SOC) relaxation for unit commitment problem are proposed in \citep{bai2015decomposition,liu2018global}. While the SOC relaxation convexifies the continuous network laws, commitment decisions need to be made, resulting in a computationally challenging in mixed-integer SOC Problem \citep{bonami2011algorithms} due to the high cost of solving conic subproblems within branch-and-bound framework. To mitigate these issues, some authors have proposed strengthening the SOCP relaxation using valid inequalities and bound tightening \citep{kocuk2016strong,coffrin2015strengthening}, as well as decomposition algorithms to separate commitment decisions and AC feasibility \citep{bai2015decomposition,Gonzalo2022AC,tuncer2022misocp}.

Despite being substantially harder to solve than DC unit commitment problem formulated as a mixed-integer linear program \citep{carrion2006computationally,conejo2018power}, SOC-relaxed AC unit commitment provides commitment that are more consistent with physical network constraints. DC unit commitment relies on linearized DC power flow, ignoring reactive power, voltage magnitudes, and losses \citep{wood2013power}, and may therefore yield commitment schedules that are infeasible or require significant corrective actions when evaluated under full AC laws \citep{bai2015decomposition,molzahn2019survey}. In contrast, SOC-based AC unit commitment captures key AC feasibility aspects while maintaining convexity in the continuous relaxation, producing tighter lower bounds and more reliable commitment decisions even if the relaxation is not exact \citep{jabr2006radial}. As a result, SOC AC unit commitment is widely viewed as a principled compromise between the scalability of DC unit commitment and the accuracy of full nonconvex AC unit commitment \citep{tuncer2022misocp}.

Developing commitment schedules that are robust against contingencies has become increasingly important in modern power system operations. Robust network-constrained unit commitment models aim to identify commitment decisions that remain feasible and cost-effective under the worst-case scenario of uncertain parameters \citep{jiang2011robust,wang2025two,bertsimas2012adaptive,amjady2016adaptive}. \citet{an2014exploring} deployed column-and-constraint generation to solve this problem. More recently, distributionally robust optimization has been introduced to unit commitment formulations, offering protection against ambiguity in probability distributions. Applications span both linear~\citep{zhao2017distributionally,duan2017data} and nonlinear settings \citep{dehghan2021distributionally}. More specifically, stochastic unit commitment with DC-based formulations have been extended to incorporate weather-driven contingencies: for example, DC-based formulations that consider probabilistic transmission line disconnections based on weather forecasts \citep{sang2019integrated,mohammadi2020tractable}. Additionally, data-driven methods have emerged, with machine learning being applied to improve both scenario generation and solution quality~\citep{mohammadi2021machine}. In addition to traditional robust and stochastic formulations, resilience can be integrated into the unit commitment framework to address the challenges posed by extreme weather events. \citep{zhao2020resilient,yang2024resilient}

Security-constrained unit commitment commonly incorporates the $N-1$ security criterion, requiring the system to remain feasible following the outage of any single critical component, such as a generator or transmission line \citep{sundar2019chance}. Extensions to the more general $N-k$ criterion, which considers simultaneous outages of up to $k$ components, have been studied to address extreme events and correlated failures, but they significantly increase computational complexity due to the combinatorial growth in the number of contingencies \citep{street2010contingency,bienstock2014chance}. As a result, practical implementations often rely on screening, contingency selection, or cutting-plane and Benders-type decomposition methods to identify only the most critical contingencies during the solution process.

In this work, we focus on robust decision that provide protection against the failure of transmission lines due to hurricanes, and the number of failure lines is not given, dislike the model for $N-k$ criterion. We emphasize that severe weather is the leading reason of power outages in the United States \citep{campbell2012weather,sang2018effective}.

\subsection{Contributions}

Considering the above literature review, the contributions of this paper are threefold:

\begin{enumerate}[label=\arabic*)]
	\item Within an adaptive robust optimization framework, we model the operation problem (multi-period optimal power flow) as a second-order-conic problem, not as a linear one. Most works reported in the literature use linear formulations. 
	\item We model hurricane trajectories (and the subsequent transmission lines disabled) to be incorporated within the proposed adaptive robust framework.
	\item To ease the computational burden of the master problem of the column-and-constraint generation algorithm, we use a novel problem-specific outer-linearization technique, which results in significant computational savings. This is the major contribution of our work. We conduct extensive computational simulations to show its relevance.
\end{enumerate}

\subsection{Paper organization}
The remainder of this paper is structured as follows. Section~\ref{sec:formulation} introduces the AC network–constrained unit commitment problem, presents its second-order cone reformulation, and derives the corresponding dual after fixing the binary decisions. Section~\ref{sec:tri-level} develops the proposed three-level adaptive robust optimization framework. Section~\ref{sec:algorithm} describes the column-and-constraint generation scheme together with the proposed outer-approximation strategy. Section~\ref{sec:experiments} reports comprehensive numerical experiments that validate the proposed formulation and demonstrate the performance of the solution approach. Section~\ref{sec:conclusion} concludes the paper.

\section{Formulation of the network-constrained unit commitment problem}\label{sec:formulation}

\subsection{Problem statements and model assumptions}

In this section, we detail the formulation of the considered unit commitment problem with the following features:
\begin{enumerate}
	\item Flows through lines are represented using relaxed SOC expressions.
	\item Lines can be disabled due to hurricane trajectories and adopt a mixed-integer linear formulation similar to that proposed by \citet{Burak2017New} regarding line switching.
	\item { Our formulation includes unserved active and reactive power (slack variables) to guarantee feasibility of the relaxed formulation, and slack variables are penalized in the objective.}
\end{enumerate}

\subsection{Primal problem: network-constrained unit commitment}

We present the AC network-constrained unit commitment problem with second-order-conic relaxed expression in Problem~(\ref{eq:acncuc}). A detailed explanation to the each constraints is provided after the formulation. For some constraints, we indicate the corresponding dual variable in parentheses after the constraint for the subsequent dual problem.

\begin{subequations}\label{eq:acncuc}
	\begin{align}
		\min\ &\sum_{t\in\mathcal{T}}\sum_{g\in\mathcal{G}} C_{g}^\mathrm{F}u_{g,t}+C_{g}^\mathrm{SU}y_{g,t}+C_{g}^\mathrm{SD}z_{g,t}
		+\sum_{t\in\mathcal{T}}\left[\sum_{g\in\mathcal{G}}C_g^\mathrm{V} p_{g,t}+\sum_{n\in\mathcal{N}} C^\mathrm{U}(p_{n,t}^\mathrm{U} + q_{n,t}^\mathrm{U})\right]\label{eq:p_obj}\\
		\mbox{s.t.}\ &y_{g,t}-z_{g,t} = u_{g,t} - u_{g,t-1},\quad \forall t=2,...,T, \forall g\in\mathcal{G} \label{eq:p_logic1}\\
		&y_{g,1}-z_{g,1} = u_{g,1} - u_{g,0},\quad \forall g\in\mathcal{G}\label{eq:p_logic1_0}\\
		&y_{g,t}+z_{g,t} \leq 1,\quad \forall t\in\mathcal{T}, \forall g\in\mathcal{G}\\
		&u_{g,t},y_{g,t},z_{g,t}\in\{0,1\},\quad \forall t=1,...,T, \forall g\in\mathcal{G}\label{eq:p_logic2}\\
		&\sum_{t=1}^{L_g} (1-u_{g,t}) = 0,\quad \forall g\in \mathcal{G}^\mathrm{T}\label{eq:p_updown1}\\
		&\sum_{t = \tilde{t}}^{\tilde{t}+T^\mathrm{U}_g-1}u_{g,t} \geq T_g^\mathrm{U} y_{g,\tilde{t}},\quad \forall \tilde{t}=L_g+1,...,T - T^\mathrm{U}_g+1,\forall g\in\mathcal{G}^\mathrm{T}\label{eq:p_updown4}\\
		&\sum_{t = \tilde{t}}^\mathrm{T}(u_{g,t} - y_{g,t}) \geq 0,\quad \forall \tilde{t}=T - T^\mathrm{U}_g+2,...,T,\forall g\in\mathcal{G}^\mathrm{T}\label{eq:p_updown6}\\
		&\sum_{t=1}^{F_g} u_{g,t} = 0,\quad \forall g\in\mathcal{G}^\mathrm{T}\label{eq:p_updown3}\\
		&\sum_{t = \tilde{t}}^{t_i+T^\mathrm{D}_g-1}(1-u_{g,t}) \geq T_g^\mathrm{D} z_{g,\tilde{t}},\quad \forall \tilde{t}=F_g+1,...,T - T^\mathrm{D}_g+1,\forall g\in\mathcal{G}^\mathrm{T}\label{eq:p_updown5}\\
		&\sum_{t = \tilde{t}}^\mathrm{T}(1 - u_{g,t} - z_{g,t}) \geq 0,\quad \forall \tilde{t}=T - T^\mathrm{D}_g+2,...,T, \forall g\in\mathcal{G}^\mathrm{T}\label{eq:p_updown2}\\
		&p_g^{\min}u_{g,t}\leq p_{g,t}\leq p_g^{\max}u_{g,t},\quad \forall t\in\mathcal{T}, g\in\mathcal{G}\quad (\overline{\gamma_{g,t}}, \underline{\gamma_{g,t}}) \label{eq:p1}\\
		&p_{g,t} \leq \bar{p}_{g,t}\leq p_g^{\max}u_{g,t},\quad \forall t\in\mathcal{T}, g\in\mathcal{G}\quad (\overline{\bar{\gamma}_{g,t}}, \underline{\bar{\gamma}_{g,t}})\\
		&q_g^{\min}u_{g,t}\leq q_{g,t} \leq q_g^{\max}u_{g,t},\quad \forall t\in\mathcal{T}, g\in\mathcal{G}\quad (\overline{\kappa_{g,t}}, \underline{\kappa_{g,t}})\label{eq:p_bd2}\\
		&p_{g,t}-p_{g,t-1} \leq R_g^\mathrm{U}u_{g,t-1}+R^\mathrm{SU}_g y_{g,t},\quad \forall t\in\mathcal{T}, g\in\mathcal{G}^\mathrm{T}\quad (\overline{\lambda_{g,t}})\label{eq:p_ramp1}\\
		&\bar{p}_{g,t} \leq p_{g,t-1} + R_g^\mathrm{U}u_{g,t-1}+R^\mathrm{SU}_g y_{g,t},\quad  t\in\mathcal{T}, g\in\mathcal{G}^\mathrm{T}\quad (\overline{\rho_{g,t}})\label{eq:p_ramp3}\\
		&p_{g,t-1}-p_{g,t} \leq R_g^\mathrm{D}u_{g,t-1}+R^\mathrm{SD}_g z_{g,t},\quad \forall t\in\mathcal{T}, g\in\mathcal{G}^\mathrm{T}\quad (\underline{\lambda_{g,t}})\label{eq:p_ramp4}\\
		&\bar{p}_{g,t} \leq p^{\max}_g(u_{g,t}-z_{g,t+1}) + z_{g,t+1}R^\mathrm{SD}_{g},\quad \forall t\in\mathcal{T}, g\in\mathcal{G}^\mathrm{T}\quad (\underline{\rho_{g,t}})\label{eq:p_ramp2}\\
		&\sum_{g\in \Omega_r}\bar{p}_{g,t} - p_{g,t} \geq R_{r,t}^\mathrm{D},\quad \forall g\in\mathcal{G}^\mathrm{T}, \forall r\in \mathcal{R}\quad (\delta_{r,t})\label{eq:p_reserve}\\
		&\sum_{g\in \Lambda_n^\mathrm{G}}p_{g,t} - p^\mathrm{D}_{n,t} + p^\mathrm{U}_{n,t} =\sum_{m\in\Lambda_n}p_{n,m,t},\quad \forall t\in\mathcal{T},n\in \mathcal{N}\quad (\mu_n,t)\label{eq:p_bal1}\\
		&\sum_{g\in \Lambda_n^\mathrm{G}}q_{g,t} - q^\mathrm{D}_{n,t} + q^\mathrm{U}_{n,t} =\sum_{m\in\Lambda_n}q_{n,m,t},\quad \forall t\in\mathcal{T}, n\in \mathcal{N}\quad (\nu_{n,t})\label{eq:p_bal2}\\
		&p_{n,m,t} = -G_{n,m}c_{n,n,t}^{m}+G_{n,m}c_{n,m,t}-B_{n,m}s_{n,m,t},\quad \forall t\in\mathcal{T}, \forall (n,m)\in\mathcal{E}\quad (\xi_{n,m, t})\label{eq:p_flow1}\\
		&q_{n,m,t} = (B_{n,m}-b^{\mathrm{shunt}}_{n,m})c_{n,n,t}^m -G_{n,m}s_{n,m,t} -B_{n,m}c_{n,m,t},\quad \forall t\in\mathcal{T}, \forall (n,m)\in\mathcal{E} \quad (\psi_{n,m, t})\label{eq:p_flow2}\\
		&f_{n,m,t} = S_{n,m, t},\quad\forall t\in\mathcal{T}, \forall (n,m)\in\mathcal{E}\quad (\zeta_{n,m, t})\label{eq:p_cap_aux}\\
		&p_{n,m,t}^2+ q_{n,m,t}^2 \leq f_{n,m,t}^2,\quad\forall t\in\mathcal{T}, \forall (n,m)\in\mathcal{E} \label{eq:p_cap}\\
		&c_{n,m,t} = c_{m,n,t}, s_{n,m,t} = -s_{m,n,t},\quad \forall t\in\mathcal{T}, \forall (n,m)\in\mathcal{E}\quad (\beta^c_{n,m,t},\beta^s_{n,m,t}) \label{eq:p_soc1}\\
		&D_{n,m,t}^{(1)} = 2c_{n,m,t}, D_{n,m,t}^{(2)} = 2s_{n,m,t},\quad \forall t\in\mathcal{T}, \forall (n,m)\in\mathcal{E}\quad (\alpha^{(k)}_{n,m,t})\label{eq:p_soc3}\\
		&D_{n,m,t}^{(3)} = c_{n,n,t}^{m} - c_{m,m,t}^{n}, D_{n,m,t}^{(4)} = c_{n,n,t}^{m} + c_{m,m,t}^{n},\quad \forall t\in\mathcal{T}, \forall (n,m)\in\mathcal{E}\label{eq:p_soc4}\\
		&(D_{n,m,t}^{(1)})^2 + (D_{n,m,t}^{(2)})^2 + (D_{n,m,t}^{(3)})^2 \leq (D_{n,m,t}^{(4)})^2,\quad \forall t\in\mathcal{T}, \forall (n,m)\in\mathcal{E} \label{eq:p_soc2}\\
		&0\leq p^\mathrm{U}_{n,t}\leq p^\mathrm{D}_{n,t},\quad \forall n\in\mathcal{N},t\in\mathcal{T}\quad (\overline{\tau_{n,t}^{p}}, \underline{\tau_{n,t}^{p}})\label{eq:p_unser1}\\
		&0\leq q^\mathrm{U}_{n,t}\leq q^\mathrm{D}_{n,t},\quad \forall n\in\mathcal{N},t\in\mathcal{T}\quad (\overline{\tau_{n,t}^{q}}, \underline{\tau_{n,t}^{q}})\label{eq:p_unser2}\\
		&-a_{n,m}V_n^{\max}V_{m}^{\max}\leq c_{n,m, t}\leq a_{n,m}V_n^{\max}V_{m}^{\max},\quad \forall t\in\mathcal{T}, \forall (n,m)\in\mathcal{E}\quad (\overline{\eta_{n,m,t}^{c}}, \underline{\eta_{n,m,t}^{c}})\label{eq:p_lin1}\\
		&-a_{n,m}V_n^{\max}V_{m}^{\max}\leq s_{n,m, t}\leq a_{n,m}V_n^{\max}V_{m}^{\max},\quad \forall t\in\mathcal{T}, \forall (n,m)\in\mathcal{E}\quad (\overline{\eta_{n,m,t}^{s}}, \underline{\eta_{nm,t}^{s}}) \label{eq:p_lin3}\\
		&a_{n,m}{V_n^{\min}}^2\leq c_{n,n,t}^m\leq a_{n,m}{V_n^{\max}}^2,\quad \forall t\in\mathcal{T}, \forall (n,m)\in\mathcal{E}\quad (\overline{\eta_{n,m,t}}, \underline{\eta_{n,m,t}})\label{eq:p_lin4}\\
		&c_{n,n, t} - {V_n^{\max}}^2(1-a_{n,m})\leq c_{n,n,t}^m,\quad \forall t\in\mathcal{T}, \forall (n,m)\in\mathcal{E}\quad (\underline{\sigma_{n,m,t}})\label{eq:p_lin5}\\
		&c_{n,n,t}^m \leq c_{n,n, t} - {V_n^{\min}}^2(1-a_{n,m}),\quad \forall t\in\mathcal{T}, \forall (n,m)\in\mathcal{E}\quad (\overline{\sigma_{n,m,t}})\label{eq:p_lin2}\\
		&{V_n^{\min}}^2\leq c_{n,n, t} \leq {V_n^{\max}}^2,\quad \forall n\in\mathcal{N},t\in\mathcal{T}\quad (\overline{\chi_{n,t}}, \underline{\chi_{n,t}})\label{eq:p2}
	\end{align}
\end{subequations}

The objective function~(\ref{eq:p_obj}) contains the startup, shutdown, and operation costs of all the units, the cost of the unserved energy, and the penalty for reactive power shedding.
{
We note that, since each reactive power load is linked to an active power load through a fixed power factor, any amount of unserved active power load directly implies a corresponding amount of unserved reactive power load; i.e., $\displaystyle q^\mathrm{U}_{n,t} = p^\mathrm{U}_{n,t}\frac{q^\mathrm{D}_{n,t}}{p^\mathrm{D}_{n,t}}$.
}
	
Constraints (\ref{eq:p_logic1})-(\ref{eq:p_logic2}) are logic constraints for units' startup, shutdown, and on/off status variables. Specifically, the constraint (\ref{eq:p_logic1}) enforces a startup ($y_{g,t}=1$) if the generating unit transitions from off to on, and a shutdown ($z_{g,t}=1$) if the generating unit transitions from on to off. The constraint (\ref{eq:p_logic1_0}) is for the first time period, where $u_{g,0}$ is a constant. The constraint (\ref{eq:p_logic2}) ensures a unit does not start up and shut down in the same time.

Minimum up and down times of thermal units are enforced by constraints (\ref{eq:p_updown1})-(\ref{eq:p_updown2}). We note that if a thermal generating unit $g \in\mathcal{G}^{\mathrm{T}}$ starts up or shuts down, it needs to be keep on/off during the minimum up/down time $T^{\mathrm{U}}_g/T^{\mathrm{D}}_g$. Additionally, in the initial status $t = 0$, the unit $g$ is required to be on/off during the first $L_g/F_g$ periods. Also, $L_g \cdot F_g = 0$, which prevents a unit to be up and down at the same time. Constraints (\ref{eq:p_updown1})/(\ref{eq:p_updown3}) enforces that unit $g$ remains on/off during the initial up/down time. Constraints (\ref{eq:p_updown4})/(\ref{eq:p_updown5}) ensures that unit $g$ remains on/off during the minimum up/down time after a startup/shutdown. When the required up/down time exceeds the total planning time range, constraints (\ref{eq:p_updown6})/(\ref{eq:p_updown2}) enforces that unit $g$ remains on/off until the last period $T$.

The actual power output of each unit has a lower bound as its minimum power output and a upper bound as its maximum available power by constraints (\ref{eq:p1})-(\ref{eq:p_bd2}). We also include the maximum available power $\bar{p}_{g,t}$ of unit $g$ in period $t$, which depends on the ramping limit and is needed to meet the reserve requirements per area, and that is upper bounded by the active power capacity of the unit.

Constraints (\ref{eq:p_ramp1})-(\ref{eq:p_ramp2}) impose to both the generation level $p_{g,t}$ and its upper bound $\bar{p}_{g,t}$ ramping limits. Specifically, constraints (\ref{eq:p_ramp1})/(\ref{eq:p_ramp3}) is the up and startup ramping limit for $p_{g,t}/\bar{p}_{g,t}$, while constraints (\ref{eq:p_ramp4})/(\ref{eq:p_ramp2}) are for the down and shutdown ramping limit.

For each reserve area $r\in\mathcal{R}$, we enforce that the maximum available active power output $\bar{p}_{g,t}$ and the active power produced $p_{g,t}$ satisfy some pre-specified reserve requirements, which are set by constraints (\ref{eq:p_reserve}).

The active and reactive power balances involving power generations, demands, unserved demands, and power flows are ensured by constraints (\ref{eq:p_bal1}) and (\ref{eq:p_bal2}), respectively. 
{
Regarding the reactive power flow balance~(\ref{eq:p_bal2}), we note that, since reactive power ``does not travel,'' adequate local voltage support throughout the power system is essential. We assume that reactive power compensators, in the form of capacitor banks and reactors, are available across the system. For tractability, these compensators are modeled as continuous resources.
}
Active and reactive power flows are expressed by (\ref{eq:p_flow1}) and (\ref{eq:p_flow2}), respectively, where we integrate the rectangular formulation \citep{jabr2006radial} and the linearization for line's switching on/off. More details are provided in Sec.~\ref{sec:express_details} including constraints (\ref{eq:p_soc1})-(\ref{eq:p_soc2}) and (\ref{eq:p_lin1})-(\ref{eq:p2}).

The transmission capacity limits, constraints (\ref{eq:p_cap}), are quadratic constraints. Here we use an auxiliary variable $f_{n,m,t}$, which is always equal to a constant, the line's capacity $S_{n,m,t}$, as (\ref{eq:p_cap_aux}), to make (\ref{eq:p_cap}) second-order-conic.

Slack variables $p_{nt}^U,q_{nt}^U$ are bounded by constraints (\ref{eq:p_unser1}) and (\ref{eq:p_unser2}).

\subsubsection{Rectangular formulation and line on/off status} \label{sec:express_details}

Constraints (\ref{eq:p_flow1}), (\ref{eq:p_flow2}), and (\ref{eq:p_soc1})-(\ref{eq:p_soc1}) are the rectangular formulation with second-order-conic relaxation \citep{jabr2006radial}. The details for the formulation can be found in, e.g., \citep{kocuk2016strong,Burak2017New,Gonzalo2022AC,dan2}. Here we provide a simple summary.

The complex voltage for a node $n$ at period $t$ is $V_{n,t} := e_{n,t} + \mathrm{i} f_{n,t}$, where $e_{n,t}$ and $f_{n,t}$ are the real and imaginary components, and $|V_n| = \sqrt{e_{n,t}^2 + f_{n,t}^2}$. For a line $(n,m)$, let $c_{n,n,t} := e_{n,t}^2  + f_{n,t}^2$, $c_{n,m,t} := e_{n,t}e_{m,t}  + f_{n,t}f_{m,t}$, and $s_{n,m,t} := e_{n,t}f_{m,t}  - e_{m,t}f_{n,t}$. Then, considering \citep[(1), (4d),(4e)]{Burak2017New}, we have (\ref{eq:p_flow1}) and (\ref{eq:p_flow2}) as the active and reactive power flow based on $c_{n,n,t}, c_{n,m,t}, s_{n,m,t}$. Also, constraints (\ref{eq:p_soc1}) and (\ref{eq:p2}) are properties for $c_{n,n,t}, c_{n,m,t}, s_{n,m,t}$, see \citep[(4f),(4g)]{Burak2017New}.

We also note that $c_{n,m,t}^2 + s_{n,m,t}^2 = c_{n,n,t}c_{m,m,t}$ is a nonconvex quadratic constraint. We relax this constraint as $c_{n,m,t}^2 + s_{n,m,t}^2 \leq c_{n,n,t}c_{m,m,t}$, which can be rewritten as a second-order-conic constraint:
\[4c_{n,m,t}^2 + 4s_{n,m,t}^2 + (c_{n,n,t} - c_{m,m,t})^2 \leq (c_{n,n,t} + c_{m,m,t})^2,\]
where four components correspond to auxiliary variables in constraints (\ref{eq:p_soc3}) - (\ref{eq:p_soc4}).

Additionally, we consider the on/off status of a line with the variable $a_{n,m} \in\{0,1\}$; hence if a line is off ($a_{n,m} = 0$), then $p_{n,m,t} = 0$ and $q_{n,m,t} = 0$. To do so, we impose constraints (\ref{eq:p_lin1}), (\ref{eq:p_lin3}), (\ref{eq:p_lin4}), in which, if $a_{n,m} = 0$, we force $c_{n,n,t}^m = 0, c_{n,m,t} =0$, and $s_{n,m,t} = 0$, where $c_{n,n,t}^m$ is an auxiliary variable for $p_{n,m,t}$ and $q_{n,m,t}$. For all $(n,m)\in\mathcal{E}$ such that $a_{n,m} = 1$, we need $c_{n,n,t}^m = c_{n,n,t}$, which is enforced by (\ref{eq:p_lin5}) as inspired by McCormick inequalities \citep{mccormick1976computability}.

\subsection{Dual problem}

{
We note that, for the sake of clarity, we first formulate a network-constrained convexified AC unit commitment problem, given as problem \eqref{eq:acncuc}. We then fix the binary variables to given values, resulting in a multi-period convexified optimal power flow problem. Finally, we derive the dual of this multi-period convexified optimal power flow problem. This is further clarified below.
}

After fixing all commitment binary variables, i.e., all $u_{g,t},y_{g,t},z_{g,t}$, and all lines' on/off status $a_{nm}$, Problem~(\ref{eq:acncuc}) turns to be a second-order-conic problem, which can be write as 
\begin{align}\label{eq:cuc}
	\arg\min\{\sum_{t\in\mathcal{T}}(\sum_{g\in\mathcal{G}}C_g^\mathrm{V} p_{gt}+\sum_{n\in\mathcal{N}} C^\mathrm{U}(p_{nt}^\mathrm{U} + q_{nt}^\mathrm{U}))\ |\ (\ref{eq:p1})-(\ref{eq:p2})\},
\end{align}
which is convex.

Except constraints (\ref{eq:p_cap}), (\ref{eq:p_soc2}), all other constraints are linear, whose duals are also linear functions. Constraints (\ref{eq:p_cap}) are (\ref{eq:p_soc2}) are second-order conic, and their conic dualities are their self-duals. For example,
\[p_{n,m,t}^2+ q_{n,m,t}^2 \leq f_{n,m,t}^2,\]
where $p_{n,m,t}, q_{n,m,t}$ and $f_{n,m,t}$ are defined in (\ref{eq:p_flow1}), (\ref{eq:p_flow2}), and (\ref{eq:p_cap_aux}), and their dual variables are $\xi_{n,m,t}, \psi_{n,m,t}$, and $\zeta_{n,m,t}$ separately. Therefore, the dual of the second-order cone is
\[\xi_{n,m,t}^2 + \psi_{n,m,t}^2 \leq \zeta_{n,m,t}^2.\]

Denominating $V_{nm}:=V_n^{\max}V_{m}^{\max}$, the dual problem of Problem~(\ref{eq:cuc}) is
\begin{subequations}\label{eq:soc_dual}
	\begin{align}
		\max\quad &\sum_{t\in\mathcal{T}}\sum_{g\in\mathcal{G}} (\underline{\gamma_{g,t}}p_g^{\min}u_{g,t}-\overline{\gamma_{g,t}}p_g^{\max}u_{g,t} - \overline{\bar{\gamma}_{g,t}}p_g^{\max}u_{g,t}+\notag \\
		&\underline{\kappa_{g,t}}q_g^{\min}u_{g,t}-\overline{\kappa_{g,t}}q_g^{\max}u_{g,t} - (R_g^\mathrm{U}u_{g,t-1}+R_g^\mathrm{SU}y_{g,t})(\overline{\lambda_{g,t}} + \overline{\rho_{g,t}}) -\notag\\
		&(R_g^\mathrm{D}u_{g,t} + R_g^\mathrm{SD}z_{g,t})\underline{\lambda_{g,t}} - p_g^{\max}(u_{g,t} -z_{g,t+1}+ z_{g,t+1}R_g^\mathrm{SD})\underline{\rho_{g,t}}) +\notag \\
		&\sum_{t\in\mathcal{T}}(\sum_{r\in \mathcal{R}}R^\mathrm{D}_t\delta_{r,t}+\sum_{n\in\mathcal{N}} (p_{n,t}^\mathrm{D}(\mu_{n,t} - \overline{\tau_{n,t}^p}) + q_{n,t}^\mathrm{D}(\nu_{n,t} - \overline{\tau_{n,t}^q}))) - \notag \\
		&\sum_{t\in\mathcal{T}}\sum_{(n,m)\in\mathcal{L}}(S_{n,m}\zeta_{n,m,t} - a_{n,m}V_{n,m}(\overline{\eta_{n,m,t}^x} + \overline{\eta_{n,m,t}^s} + \underline{\eta_{n,m,t}^x} + \underline{\eta_{n,m,t}^s}) + \notag \\
		&a_{n,m}{V_n^{\min}}^2\underline{\eta_{n,m,t}} - a_{n,m}{V_n^{\max}}^2\overline{\eta_{n,m,t}} - {V_n^{\max}}^2(1-a_{n,m})\underline{\sigma_{n,m,t}} +\notag \\
		&{V_n^{\min}}^2(1-a_{n,m})\overline{\sigma_{n,m,t}}+ {V_n^{\min}}^2\underline{\chi_{n,t}}-{V_n^{\max}}^2\overline{\chi_{n,t}})\label{eq:dual_obj}\\
		\mbox{s.t.}\quad &\underline{\gamma_{g,t}} - \overline{\gamma_{g,t}} - \underline{\bar{\gamma}_{g,t}} - \overline{\lambda_{g,t}} + \overline{\lambda_{g,t+1}} + \overline{\rho_{g,t+1}}+ \underline{\lambda_{g,t}}\notag\\
		& - \underline{\lambda_{g,t+1}} - \delta_{r,t}+\mu_{n,t} = C_g^\mathrm{V},\ \forall t= 1,..., T-1,\ g\in\mathcal{G}^\mathrm{T}\ (p_{g,t}) \label{eq:dual_y1}\\
		&\underline{\gamma_{g,t}} - \overline{\gamma_{g,t}} +\mu_{n,t} = C_g^\mathrm{V},\ \forall t\in\mathcal{T}, g\in\mathcal{G}^\mathrm{R}\ (p_{g,t})\\
		&\underline{\bar{\gamma}_{g,t}} - \overline{\bar{\gamma}_{g,t}}-\overline{\rho_{g,t}}-\underline{\rho_{g,t}} + \delta_{r,t} =0,\ \forall t\in\mathcal{T}, g\in\mathcal{G}^\mathrm{T} \ (\bar{p}_{g,t})\\
		&\underline{\bar{\gamma}_{g,t}} - \overline{\bar{\gamma}_{g,t}} + \delta_{r,t} =0,\ \forall t\in\mathcal{T}, g\in\mathcal{G}^\mathrm{R} \ (\bar{p}_{g,t})\\
		&\underline{\kappa_{g,t}} - \overline{\kappa_{g,t}} + \nu_{n,t} = 0,\ \forall t\in\mathcal{T}, g\in\mathcal{G}\ (q_{g,t})\\
		&\mu_{n,t} + \underline{\tau^{p}_{n,t}} - \overline{\tau^{p}_{n,t}} = C^\mathrm{U},\ \forall t\in\mathcal{T}, n\in\mathcal{N}\ (p_{n,t}^{U})\\
		&\nu_{n,t} + \underline{\tau^{q}_{n,t}} - \overline{\tau^{q}_{n,t}} = C^\mathrm{U},\ \forall t\in\mathcal{T}, n\in\mathcal{N}\ (q_{n,t}^{U}) \label{eq:dual_y2}\\
		&\underline{\eta_{n,m,t}} - \overline{\eta_{n,m,t}} + \underline{\sigma_{n,m,t}} - \overline{\sigma_{n,m,t}} + \alpha_{n,m,t}^{(3)}-\alpha_{m,n,t}^{(3)} + \notag\\
		&\alpha_{n,m,t}^{(4)} +\alpha_{m,n,t}^{(4)}-G_{n,m}(\xi_{n,m,t} -\mu_{n,t})+\notag\\
		&(B_{n,m}-b_{n,m}^{shunt})(\psi_{n,m,t}-\nu_{n,t})=0,\ \forall t\in\mathcal{T}, (n,m)\in\mathcal{E} \ (c_{n,n,t}^{m}) \label{eq:dual_z1}\\
		&\underline{\chi_{n,t}}-\overline{\chi_{n,t}} + \sum_{m\in\Lambda_n}(\overline{\sigma_{n,m,t}}-\underline{\sigma_{n,m,t}}) = 0,\ \forall t\in\mathcal{T} \ (c_{n,n,t})\\
		&\underline{\eta_{n,m,t}^{x}}-\overline{\eta_{n,m,t}^{x}} + G_{n,m}(\xi_{n,m,t} -\mu_{n,t}) - B_{n,m}(\psi_{n,m,t}-\nu_{n,t})\notag\\
		&+ \beta_{n,m,t}^{x} -\beta_{m,n,t}^x+2\alpha_{n,m,t}^{(1)} = 0,\ \forall t\in\mathcal{T}, (n,m)\in\mathcal{E} \ (c_{n,m,t})\\
		&\underline{\eta_{n,m,t}^{s}}-\overline{\eta_{n,m,t}^{s}} - B_{n,m}(\xi_{n,m,t} -\mu_{n,t}) - G_{n,m}(\psi_{n,m,t}-\nu_{n,t})\notag\\
		&+ \beta_{n,m,t}^{s} +\beta_{m,n,t}^s+2\alpha_{n,m,t}^{(2)} = 0,\ \forall t\in\mathcal{T}, (n,m)\in\mathcal{E} \ (s_{n,m,t})  \label{eq:dual_z2}\\
		&\xi_{n,m,t}^2 + \psi_{n,m,t}^2 \leq \zeta_{n,m,t}^2, \ \forall t\in\mathcal{T}, (n,m)\in\mathcal{E} \label{eq:dual_g1}\\
		&(\alpha_{n,m,t}^{(1)})^2 + (\alpha_{n,m,t}^{(2)})^2+ (\alpha_{n,m,t}^{(3)})^2\leq (\alpha_{n,m,t}^{(4)})^2, \ \forall t\in\mathcal{T}, (n,m)\in\mathcal{E} \label{eq:dual_g2}\\
		&\overline{\gamma_{g,t}}, \underline{\gamma_{g,t}}, \overline{\bar{\gamma}_{g,t}}, \underline{\bar{\gamma}_{g,t}}, \overline{\kappa_{g,t}}, \underline{\kappa_{g,t}}, \overline{\lambda_{g,t}}, \underline{\lambda_{g,t}}, \overline{\rho_{g,t}}, \underline{\rho_{g,t}}, \delta_{g,t}, \overline{\tau_{n,t}^{p}}, \underline{\tau_{n,t}^{p}}, \overline{\tau_{n,t}^{q}}, \underline{\tau_{n,t}^{q}}\geq 0\\ &\overline{\eta_{n,m,t}^{x}}, \underline{\eta_{n,m,t}^{x}}, \overline{\eta_{n,m,t}^{s}}, \underline{\eta_{n,m,t}^{s}}, \overline{\eta_{n,m,t}}, \underline{\eta_{n,m,t}},
		\overline{\sigma_{n,m,t}}, \underline{\sigma_{n,m,t}},
		\overline{\chi_{n,t}}, \underline{\chi_{n,t}}\geq 0.
	\end{align}
\end{subequations}

\section{Tri-Level robust formulation with uncertain hurricane trajectories} \label{sec:tri-level}

The tri-level problem, which is the key component of this paper is described below.

\subsection{Compact tri-level formulation}

The adaptive robust formulation is compactly stated as

\begin{subequations}\label{eq:compact_three}
	\begin{align}
		\min_{\boldsymbol{x}}\ &c_x^T\boldsymbol{x} + \max_{a\in \Omega_{\mathrm{traj}}}\ \min_{(\boldsymbol{y,z})\in\Omega(\boldsymbol{x},a)} c_y^T\boldsymbol{y} + c_z^T\boldsymbol{z}\\
		\mbox{s.t.}\ & F\boldsymbol{x} \leq f,\ \boldsymbol{x} \mbox{ binary} \label{eq:primal_1}
	\end{align}
\end{subequations}
where 
\begin{subequations}
	\begin{align}
		\Omega(\boldsymbol{x}, a):=\{&H_y\boldsymbol{y}(a) + H_z\boldsymbol{z}(a)\leq h(a), \label{eq:omega_1}\\
		&A\boldsymbol{x} +B\boldsymbol{y}(a) + C\boldsymbol{z}(a)\leq d, \label{eq:omega_2}\\
		&\boldsymbol{g}^{cap}(\boldsymbol{z}(a), a) \leq \boldsymbol{s}, \label{eq:omega_3}\\
		&\boldsymbol{g}^{soc}(\boldsymbol{z}(a), a) \leq 0\}.
	\end{align}
\end{subequations}

The variable $\boldsymbol{x}$ contains all commitment variables, i.e., $u_{gt},y_{gt},z_{gt}$ from Problem~(\ref{eq:acncuc}), the variable $\boldsymbol{y}(\cdot)$ contains all dispatch variables, e.g., $p_{gt}, \bar{p}_{gt}, q_{gt}$, and $\boldsymbol{z}(\cdot)$ contains all network variables, e.g., $p_{nm, t}, q_{nm, t}, c_{nm,t},$ $s_{nm,t}, c_{nn,t}^m$.

Constraint (\ref{eq:primal_1}) covers constraints (\ref{eq:p_logic1}) - (\ref{eq:p_updown2}); constraint (\ref{eq:omega_2}) covers constraints (\ref{eq:p1}) - (\ref{eq:p_bal2}); and constraint (\ref{eq:omega_1}) covers constraints (\ref{eq:p_flow1}) - (\ref{eq:p2}).

There are two nonlinear constraint vectors: $\boldsymbol{g}^{cap}(\boldsymbol{z}(\cdot), \cdot)$ is a vector for line capacity limits, and constraints (\ref{eq:p_cap}), where $\boldsymbol{s}:=(S_{nm}^2)_{(n,m)\in\mathcal{E}}$; and $\boldsymbol{g}^{soc}(\boldsymbol{z}(\cdot), \cdot)$ is a vector for second-order cone relaxations, constraints (\ref{eq:p_soc2}), i.e., $c_{nm}^2 + s_{nm}^2 \leq c_{nn}^mc_{mm}^n$.

\subsection{Compact primal and dual third-level problem}

The third-level problem of model~(\ref{eq:compact_three}) can be represented as 
$$\bar{\Omega}(\boldsymbol{x}, a):=\min_{\boldsymbol{y,z}\in\Omega(\boldsymbol{x},a)}\{c_y^T\boldsymbol{y}+c_z^T\boldsymbol{z}\}.$$
Its dual can be written compactly as
\begin{subequations}
	\begin{align}
		S(\boldsymbol{x}, a):= \max\ &(A\boldsymbol{x} - d)^T\lambda - h(a)^T\sigma + \boldsymbol{s}^T\eta\\
		\mbox{s.t}\ & H_y^T\lambda + B^T\sigma = c_y, \label{eq:dual_1}\\
		&H_z^T\lambda + C^T\sigma = c_z, \label{eq:dual_2}\\
		&\boldsymbol{g}^{cap}(\eta, a) \leq \boldsymbol{s}, \label{eq:dual_3}\\
		&\boldsymbol{g}^{soc}(\sigma, a) \leq 0, \label{eq:dual_4},\\
		&\lambda,\sigma,\eta\geq 0
	\end{align}
\end{subequations}
where $\sigma, \lambda$, and $\eta$ are dual variables of constraints (\ref{eq:omega_1}), (\ref{eq:omega_2}), and (\ref{eq:omega_3}), respectively.

It is important to note that constraint (\ref{eq:dual_1}) covers constraints (\ref{eq:dual_y1})-(\ref{eq:dual_y2}); constraint (\ref{eq:dual_2}) covers constraints (\ref{eq:dual_z1}) - (\ref{eq:dual_z2}); constraint (\ref{eq:dual_3}) covers constraint (\ref{eq:dual_g1}); and constraint (\ref{eq:dual_4}) covers constraint (\ref{eq:dual_g2}).

\subsection{Hurricane selection problem: second-level problem}

The second-level problem of model~(\ref{eq:compact_three}) is $\max_{a\in\Omega_{\mathrm{traj}}}\{\bar{\Omega}(\boldsymbol{x}, a)\}$, where $a$ represents the on/off status of lines. Due to a hurricane, line $(n,m)\in\mathcal{E}$ may become disabled, i.e. $a_{n,m} = 0$. Specifically, consider there are $n_{\mathrm{traj}}$ different hurricane trajectories. For the $k$-th trajectory, there are $a_{n,m}^{(k)}\in\{0,1\}$ for all $(n,m)\in\mathcal{E}$. If $a_{n,m}^{(k)} = 0$, line $(n,m)$ is disabled; otherwise it is not, and $a_{n,m}^{(k)} = 1$. For convenience, let $a_{n,m}^{(0)} := 1$ for all $(n,m)\in\mathcal{E}$, which represents the scenario with no hurricane. Let $\Omega_{\mathrm{traj}}$ be the set including all possible $a_{n,m}$ under different trajectories,
\[\Omega_{\mathrm{traj}}:=\{a\ |\ a_{n,m}=1-\sum_{k=1}^{n_{\mathrm{traj}}}\omega_k(1-a_{n,m}^{(k)}),\sum_{k=1}^{n_{\mathrm{traj}}} \omega_k\leq 1, \omega_k\in\{0,1\}\},\]
where $\omega_k\in\{0,1\}$ is a trajectory indicator, e.g., $\omega_k = 1$ implies that the $k$-th trajectory happens; otherwise, $\omega_k = 0$. We assume that at most one trajectory realizes. Detecting the worst hurricane trajectories from $\Omega_{\mathrm{traj}}$, i.e. find the $a$ maximizing the power generation costs, is the second-level problem.

\subsection{Merged second- and third-level problems}

With the description of $\Omega_{\mathrm{traj}}$,  we can merge the second-level problem and the dual of the third-level problem as
\begin{subequations}
	\begin{align}
		R(\boldsymbol{x}):=\max\ &(A\boldsymbol{x} - d)^T\lambda - h(a)^T\sigma + \boldsymbol{s}^T\eta\\
		\mbox{s.t}\ & H_y^T\lambda + B^T\sigma = c_y,\\
		&H_z^T\lambda + C^T\sigma = c_z,\\
		&\boldsymbol{g}^{cap}(\eta, a) \leq \boldsymbol{s},\\
		&\boldsymbol{g}^{soc}(\sigma, a) \leq 0,\\
		&a\in \Omega_{\mathrm{traj}}.
	\end{align}
\end{subequations}
Since $h(a)$ is a linear function, there are bilinear terms involving binary and continuous variables in the objective function of $R(\boldsymbol{x})$: $h(a)^T\sigma, a\in \{0,1\}^{|a|}, \sigma\in \mathbb{R}^{|\sigma|}$. These bilinear terms have the form $h(a_{nm})\sigma_{nm, t}$, $\sigma_{nm, t}\geq 0$. They can be exactly linearized as:
\begin{equation}\label{eq:bigM}
	\begin{aligned}
		0\leq &a_{nm}\sigma_{nm,t} \leq Ma_{nm}\\
		\sigma_{nm,t} - M(1-a_{nm}) \leq &a_{nm}\sigma_{nm,t} \leq \sigma_{nm,t}
	\end{aligned}
\end{equation}
We denote the linearized version of $R(\boldsymbol{x})$ as $\bar{R}(\boldsymbol{x})$.

\subsection{First-level Problem}

The first-level problem of the model~(\ref{eq:compact_three}) can be represented as 
$$\min_{\boldsymbol{x}} \{c_x^T\boldsymbol{x} + R(\boldsymbol{x})\ |\ F\boldsymbol{x}\leq f, \boldsymbol{x} \mbox{ binary}\}$$.

Binary variables $\boldsymbol{x}$ are commitment variables, and $F\boldsymbol{x}\leq f$ represents logic constraints and ramping limits. $R(\boldsymbol{x})$ is the recourse function from the merged second- and third-level problems.

\section{Solution method}\label{sec:algorithm}

\subsection{Column-and-Constraint generation algorithm}
In this section, we adopt the column-and-constraint generation algorithm \citep{Bo2017Solving} to solve the resulting bi-level problem. Let $\mathcal{K}$ be an index set and $\mathcal{K}\subseteq \{0,...,n_{\mathrm{traj}}\}$. The master problem related to problem~(\ref{eq:compact_three}) is:
\begin{equation}
	\begin{aligned}
		M(\mathcal{K}):= \min_{\boldsymbol{x}}\ &c_x^T\boldsymbol{x} + R\\
		\mbox{s.t.}\ &F\boldsymbol{x} \leq f,\ \boldsymbol{x} \mbox{ binary},\\
		&R\geq c_y^T\boldsymbol{y}^{(k)} + c_z^T\boldsymbol{z}^{(k)},\ \forall k\in\mathcal{K},\\
		&(\boldsymbol{y}^{(k)}, \boldsymbol{z}^{(k)})\in\Omega(\boldsymbol{x}, a^{(k)}),\  \forall k\in\mathcal{K}
	\end{aligned}
\end{equation}
The subproblem is the merged second- and third-level problem $\bar{R}(\boldsymbol{x})$. The column-and-constraint generation algorithm is provided as Alg.~\ref{Alg:ccg}.

\begin{algorithm}[!htbp]
	\SetAlgoLined
	\LinesNumbered
	\SetKwRepeat{Do}{do}{while}
	\SetKwInput{Input}{Input}
	\SetKwInput{Output}{Output}
	\Input{Problem~(\ref{eq:compact_three}), $\mbox{max\_iter}$, tolerance $\epsilon$.}
	\Output{$\boldsymbol{x}$}
	Set $\mathrm{LB} := M(\{0\})$, and $\boldsymbol{x}$ be the solution\;
	Set $\mathrm{UB} := \bar{R}(\boldsymbol{x})+c_x^T\boldsymbol{x}$, and detect $\omega_{(k_0)}= 1$ from the solution\;
	Initialize $\mathcal{K}:=\emptyset$, $gap := \frac{\mathrm{UB} - \mathrm{LB}}{\mathrm{UB}}$\;
	\For{$i=1:\mathrm{max\_iter}$}{
		Set $\mathcal{K}:=\mathcal{K}\cup\{k_{i-1}\}$\;
		Set $\mathrm{LB} := M(\{\mathcal{K}\})$, and $\boldsymbol{x}$ be the solution\;
		Set $\mathrm{UB} := \min\{\mathrm{UB}, \bar{R}(\boldsymbol{x})+c_x^T\boldsymbol{x}\}$, and detect $\omega_{(k_i)} = 1$ from the solution\;
		Set $gap:= \frac{\mathrm{UB - LB}}{\mathrm{UB}}$\;
		\If{$gap < \epsilon$}{
			\textbf{Break}\;
		}
	}
	\textbf{return} $\boldsymbol{x}$
	\caption{Column-and-Constraint Generation Algorithm}
	\label{Alg:ccg}
\end{algorithm}

On line 1 of Algorithm~\ref{Alg:ccg}, we use the solution to AC-NCUC problem with no line disabled ($M(\{0\})$) as the initial solution.

On lines 2 and 7, we solve the subproblem and detect the ``worst'' trajectory $k_i$ such that $\omega_{k_i} = 1$ from $\Omega_{\mathrm{traj}}$.

On line 5, we add the last detected trajectory (line 7) to $\mathcal{K}$.

\subsection{Outer Approximation} \label{sec:oicpa}
More often than not the master problem $M(\mathcal{K})$ is challenging to solve, especially if the size of $\mathcal{K}$ is large. Therefore, we design an outer approximation method providing linear relaxation to the quadratic line capacity constraint~(\ref{eq:p_cap}) and the second-order cone constraint~(\ref{eq:p_soc2}), both convex. To solve $M(\mathcal{K})$, we adopt a technique based on linear cuts proposed by \citet{dan1}: 

\begin{enumerate}[label=\arabic*)]
	\item Regarding the line capacity constraint $p_{nm,t}^2 + q_{nm,t}^2\leq S_{nm}^2$, if a solution $(\bar{p}_{nm,t}, \bar{q}_{nm,t})$ violates this constraint, we add the linear cut
	\[\bar{p}_{nm,t}p_{nm,t} + \bar{q}_{nm,t}q_{nm,t} \leq S_{nm}\|(\bar{p}_{nm,t}, \bar{q}_{nm,t})\|_2;\]
	\item Regarding the second-order cone constraint $c_{nm,t}^2 + s_{nm,t}^2\leq c_{nn,t}^mc_{mm,t}^n$, if a solution $(\bar{c}_{nm,t}, \bar{s}_{nm,t}$,  $\bar{c}_{nn,t}^m, \bar{c}_{mm,t}^n)$ violates this constraint, we add the cut
	\[4\bar{c}_{nm,t}c_{nm,t}+4\bar{s}_{nm,t}s_{nm,t} + (\bar{c}_{nn,t}^{m} - \bar{c}_{mm,t}^{n} - n_0)c_{nn,t}^m - (\bar{c}_{nn,t}^{m} - \bar{c}_{mm,t}^{n} + n_0)c_{mm,t}^m\leq 0,\]
	where $n_0 = \|2\bar{c}_{nm,t}, 2\bar{s}_{nm,t}, \bar{c}_{nn,t}^m, \bar{c}_{mm,t}^n\|_2$. 
\end{enumerate}

{ Let $\bar{M}(\mathcal{K}, \mathcal{C})$ be $M(\mathcal{K})$ with all nonlinear constraints relaxed, and including cuts from $\mathcal{C}$. Let $\bar{\Omega}(\boldsymbol{x}, a) :=\min_{\boldsymbol{y,z}\in\Omega(\boldsymbol{x}, a)}\{c_y^T\boldsymbol{y}+c_z^T\boldsymbol{z}\}$, the outer-inner cutting-plane method is provided in Alg.~\ref{Alg:oa}. We use $\mathrm{UB}$ and $\mathrm{LB}$ to represent the upper bound and lower bound of the master problem $\mathcal{K}$. We update $\mathrm{LB}$ as the objective value function of $\bar{M}(\mathcal{K}, \mathcal{C})$, which is an outer approximation; thus $\mathrm{LB} \leq M(\mathcal{K})$. Also, we update $\mathrm{UB}$ as the minimum among objective values of $\bar{\Omega}(\boldsymbol{x}, a^{(k)})$, which is an inner approximation to $M(\mathcal{K})$ that ignores other scenarios from $\mathcal{K}$ and fixes $\boldsymbol{x,y,z}$; thus $\mathrm{UB} \geq M(\mathcal{K})$. The algorithm will terminate when the upper and lower bounds are closed enough as line 3 of Alg.~\ref{Alg:oa}. The reported solution is from the inner approximation of the master problem; thus it is guaranteed feasible.}

\begin{algorithm}[!htbp]
	\SetAlgoLined
	\LinesNumbered
	\SetKwRepeat{Do}{do}{while}
	\SetKwInput{Input}{Input}
	\SetKwInput{Output}{Output}
	\Input{$M(\mathcal{K})$, $\epsilon$, $\epsilon_{tol}$, $p_{cut}$, and $\epsilon_{par}$.}
	\Output{$\boldsymbol{x}$ and $obj$}
	Set $\mathrm{LB}$ be the lower bound by solving $\bar{M}(\mathcal{K}, \emptyset)$, and $\boldsymbol{x,y,z}$ be the solution\;
	Initialize $\mathcal{C}:=\emptyset$, $\mathrm{UB} := +\infty$\;
	\While{$\frac{\mathrm{UB}-\mathrm{LB}}{\mathrm{UB}} <\epsilon \mathrm{\ or\ no\ new\ cuts\ are\ detected}$}{
		Check for $\epsilon_{tol}$-violated $\boldsymbol{g}^{soc}(\boldsymbol{y}, a^{(k)})\leq 0,\ \forall k\in\mathcal{K}$\;
		Find $p_{cut}$ inequalities with the most violations and generate cuts for them\;
		Push the generated cut to $\mathcal{C}$ if it is not $\epsilon_{par}$-parallel to cuts in $\mathcal{C}$\;
		\For{$k\in\mathcal{K}$}{
			$\mathrm{UB} := \min\{\mathrm{UB}, \bar{\Omega}(\boldsymbol{x}, a^{(k)})+c_x^T\boldsymbol{x}\}$\;
			Check for active $\boldsymbol{g}^{cap}(\cdot, a^{(k)})\leq \boldsymbol{s}$ by solving $\bar{\Omega}(\boldsymbol{x}, a^{(k)})$\;
			Check for $\epsilon_{tol}$-violated $\boldsymbol{g}^{cap}(\boldsymbol{z}, a^{(k)})\leq \boldsymbol{s}$ only for active inequalities and generate cuts for them\;
			Push the generated cut to $\mathcal{C}$ if it is not $\epsilon_{par}$-parallel to cuts in $\mathcal{C}$\;
		}
		Set $\mathrm{LB}$ be the lower bound by solving $\bar{M}(\mathcal{K}, \mathcal{C})$, and $\boldsymbol{x,y,z}$ be the solution\;
	}
	\textbf{return} $\boldsymbol{x}, obj:=\mathrm{UB}$\;
	\caption{Outer-Inner Cutting-Plane Method}
	\label{Alg:oa}
\end{algorithm}

Alg.~\ref{Alg:oa} works as follows:

On line 4, we detect violating inequalities of that form $g_{i}^{soc}(\boldsymbol{y}, a^{(k)})\geq \epsilon_{tol}$.

On lines 5 and 6, we select the inequalities with the top $p_{cut}$ percentage violation and compute cuts.

We note that two cuts, i.e., $c_1^Tx\leq 0$ and $c_2^Tx\leq 0$, are not $\epsilon_{par}$-parallel when the cosine of the angle formed by $c_1/\|c_1\|$ and $c_2/\|c_2\|$ is less than or equal to $1-\epsilon_{par}$.

On line 7, we add new cuts if they are not $\epsilon_{par}$-parallel. This step prevents adding too many ``too parallel'' cuts, which can slow down the progress to the solution.

On line 9, extensive computational experiments reveal that $\bar{\Omega}(\boldsymbol{x}, a^{(k)})$, which is a second-order-conic program with only one scenario, can be solved efficiently.

On lines 10 and 11, we adopt an active-set strategy that only considers binding inequalities based on the solution of the inner problem. Such an active-set strategy can reduce the number of added cuts and reduce the solution time on line 15.

We note that the proposed method (Algorithm~\ref{Alg:oa}) is different than that of \citep{dan1} in the following ways:

\begin{enumerate}[label=\arabic*)]
	\item  Alg.~\ref{Alg:oa} leverages both outer and inner approximations and stops when the solutions of the two approximations are close enough, while \citep[Algorithm 1]{dan1} stops when the objective value does not change in five continuous iterations.
	\item  Alg.~\ref{Alg:oa} adopts an active set strategy in lines 10 and 11 to reduce the size of the linearizations, which does not impact the converge guarantee. However, \citet[Algorithm 1]{dan1} eliminates inactive older cuts, which may cause divergence.
	\item Alg.~\ref{Alg:oa} is set for mixed-integer second-order-conic problems, while \citep[Algorithm 1]{dan1} is designed for continuous second-order-conic programs.
	\item Alg.~\ref{Alg:oa} returns a solution from the inner approximation, which has a feasibility guarantee, and the feasibility tolerance is the solver's tolerance, while the solution from \citep[Algorithm 1]{dan1} pertain to the outer approximation, whose feasibility is not guaranteed.
\end{enumerate}

As we mentioned in the point 4) from the comparison with \citep[Algorithm 1]{dan1} above, Alg.~\ref{Alg:oa} returns a feasible solution. Additionally, the objective value of the inner approximation and the lower bound of outer approximation provide valid upper bound and lower bound, separately, for the original problem $M(\mathcal{K})$; thus, once Alg.~\ref{Alg:oa} converges, it returns a feasible solution $M(\mathcal{K})$ with an optimality gap smaller than $\epsilon$.

\section{Numerical experiments}\label{sec:experiments}

{
\textbf{Data} In this section, we numerically validate the proposed solution approach using both the IEEE 24-bus system, and the Central Illinois 200-bus test system \citep{birchfield2016grid}. All data can be found in \citep{flores2022scheduling}, Appendix E. We note that we corrected a mistake in the data of IEEE 24-bus system regarding that susceptance of line $(6, 10)$, and made it as line $(7, 8)$. To guarantee reactive power feasibility under hurricane conditions, we set ``reactive power generators'' for nodes without any generator, which can either consume or produce reactive power as needed.

\textbf{Hardware \& Software} We conduct all experiments on an AMD Ryzen Windows64 machine including an R5-5600G 3.90GHz CPU with 6 physical cores, 12 logical processors, and 16 GB RAM. We build models and solve them with JuMP v1.25.0 \citep{Lubin2023} in Julia v1.11.5 \citep{Julia-2017}. We use two commercial solvers: Gurobi 12.0.1 \citep{gurobi} and Mosek 11.0.9 \citep{mosek} with their default setting. We solve all continuous second-order cone programs with Mosek for a better numerical robustness and solve all mixed-integer programs with Gurobi for a faster solution. We also use an open-source solver IPOPT 3.14.17 \citep{wachter2006implementation} to locally solve nonconvex programs.

\textbf{Parameters} We consider a number of time period $T=24$, i.e. the hours of a day. For the value of $M$ in Equation~(\ref{eq:bigM}), we first solve the original UC~(\ref{eq:acncuc}) without line off to get commitments, and then solve the dual problem~(\ref{eq:soc_dual}) to get dual solutions $\sigma^*$. Set $M = 100\times \max\{\sigma^*\}$. The convergence tolerance is set as 
$\epsilon := 10^{-4}$
 .}

\subsection{Case study: IEEE 24-bus power system}

\begin{figure}[!htbp]
	\centering
	\includegraphics[width=0.5\linewidth]{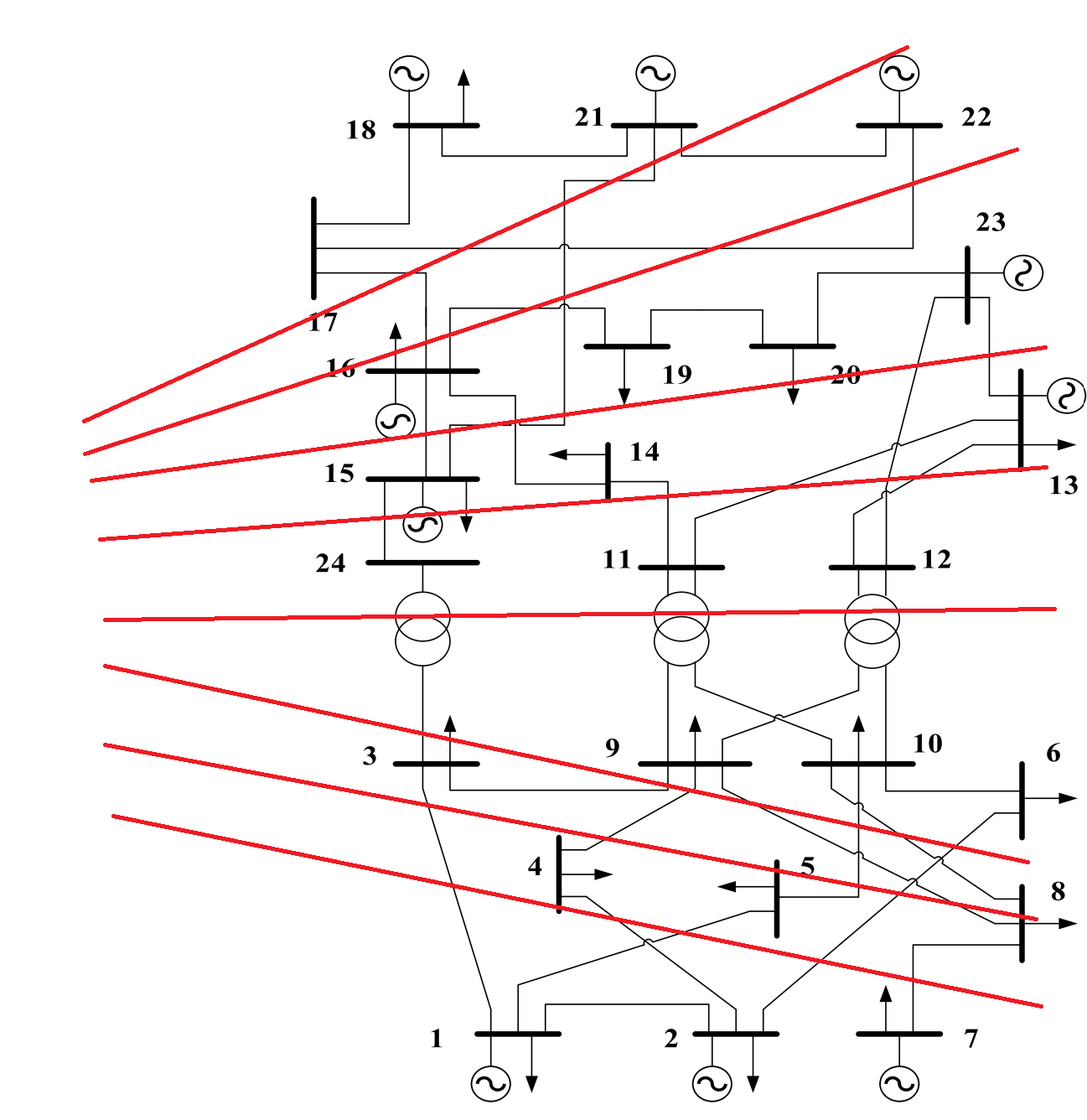}
	\caption{IEEE 24-bus system with hurricane trajectories in red}
	\label{fig:24bus}
\end{figure}

We illustrate the considered hurricane trajectories in Fig.~\ref{fig:24bus}, where each red line is a trajectory that disconnects all through lines and splits the network into two independent parts. Details of trajectories can be found in \ref{app:simulation}

Considering Alg.~\ref{Alg:ccg}, we first get an initial solution with all lines in place using Gurobi and report results in Iter. 0 of Tab.~\ref{tab:24bus}. Lower and upper bounds $\mathrm{LB}$ and $\mathrm{UB}$ are obtained by solving the master problem and subproblem separately as stated in Alg.~\ref{Alg:ccg}, and M. Runtime and S. Runtime are the corresponding runtimes for the master problem solution and subproblem solution, respectively. Additionally, as the number of iteration increases, the size of the master problem grows, which increases the solution time; thus we also present the number of nonzero elements for different master problems in different iterations as M. NNZ in Tab.~\ref{tab:24bus}.

\begin{table}[!htbp]
	\centering
	\setlength{\tabcolsep}{1.7mm}
	\begin{threeparttable}
\renewcommand{\arraystretch}{1.2}
		\begin{tabular}{l|ccc|cc|c}
			\toprule[1pt]
			Iteration & LB &M. Runtime (s) &M. NNZ &UB &S. Runtime (s) & Gap\\
			\hline
			Iter. 0 &$9.9181\times 10^{5}$ &24 &80631   &$2.2653\times 10^{8}$ &29 &99.56\%\\
			Iter. 1 &$1.8452\times 10^{8}$ &51 &98296   &$2.2653\times 10^{8}$ &45 &18.55\%\\
			Iter. 2 &$1.8474\times 10^{8}$ &109 &183497 &$2.2654\times 10^{8}$ &24 &18.45\%\\
			Iter. 3 &$1.8485\times 10^{8}$ &160 &268702 &$1.8485\times 10^{8}$ &22 &0.0\%\footnote{2}\\
			\bottomrule[1pt]
		\end{tabular}
		\begin{tablenotes}
			\footnotesize
			\item[1] The zero gap of Iter. 2 is from rounding the gap to 4 decimal places.
		\end{tablenotes}
	\end{threeparttable}
	\caption{Alg.~\ref{Alg:ccg} performance on IEEE 24-bus power system}
	\label{tab:24bus}
\end{table}

We report commitments from the initial (not robust) model and the robust model, and detect the scenario of the worst hurricane. Then we compare the performance of the two groups of commitments on the corresponding worst scenario in Table~{\ref{tab:24_comp}}. Commitment + Served Energy Cost is the cost of commitment and active power generation, while Unserved Energy Cost is the cost of the unserved active power and the penalty for reactive power shedding.

\begin{table}[!htbp]
	\centering
	\setlength{\tabcolsep}{1.7mm}
	\begin{threeparttable}
\renewcommand{\arraystretch}{1.5}
		\begin{tabular}{l|cc}
			\toprule[1pt]
			Commitment &Commitment + Served Energy Cost  &Unserved Energy Cost\\
			\hline
			Initial (Not Robust) Model & $7.9401\times 10^{5}$ &$2.2574\times 10^{8}$\\
			Robust Model               & $1.0408\times 10^{6}$ &$1.8381\times 10^{8}$\\
			\bottomrule[1pt]
		\end{tabular}
	\end{threeparttable}
	\caption{Initial (no robust) SCUC vs. Robust SCUC: Cost breakdown on IEEE 24-bus system}
	\label{tab:24_comp}
\end{table}

\begin{figure}[htbp]
	\centering
	\begin{tikzpicture}
		\begin{axis}[
			xlabel={Generating unit},
			ylabel={Time (h)},
			xtick={1,3,5,7,9,11,13,15,17,19,21,23,25,27,29,31,33},
			ytick={1,2,3,4,5,6,7,8,9,10,11,12,13,14,15,16,17,18,19,20,21,22,23,24},
			xmin=0.5, xmax=33.5,
			ymin=0.5, ymax=24.5,
			y dir=reverse, 
			enlargelimits=false,
			axis on top,
			width=14cm,
			height=10cm,
			grid=none,
			legend style={at={(0.5,1.05)}, anchor=south, legend columns=-1, column sep = 8pt},
			tick label style={font=\small},
			label style={font=\small}
			]

			\addplot[
			only marks,
			mark=square*,
			mark options={draw=blue, fill=white},
			mark size=2.5pt
			] table {
				x y
				3 1
				3 2
				3 3
				3 4
				3 5
				3 6
				3 7
				3 8
				3 9
				3 10
				3 11
				3 12
				3 13
				3 14
				3 15
				3 16
				3 17
				3 18
				3 19
				3 20
				3 21
				3 22
				3 23
				3 24
				4 1
				4 2
				4 3
				4 4
				4 5
				4 6
				4 7
				4 8
				4 9
				4 10
				4 11
				4 12
				4 13
				4 14
				4 15
				4 16
				4 17
				4 18
				4 19
				4 20
				4 21
				4 22
				4 23
				4 24
				7 1
				7 2
				7 3
				7 4
				7 5
				7 6
				7 7
				7 8
				7 9
				7 10
				7 11
				7 12
				7 13
				7 14
				7 15
				7 16
				7 17
				7 18
				7 19
				7 20
				7 21
				7 22
				7 23
				7 24
				8 6
				8 7
				8 8
				8 9
				8 10
				8 11
				8 12
				8 13
				8 14
				8 15
				8 16
				8 17
				8 18
				8 19
				8 20
				8 21
				8 22
				8 23
				8 24
				9 1
				9 2
				9 3
				9 4
				9 5
				9 6
				9 7
				9 8
				9 9
				9 10
				9 11
				9 12
				9 13
				9 14
				9 15
				9 16
				9 17
				9 18
				9 19
				9 20
				9 21
				9 22
				9 23
				9 24
				11 5
				11 6
				11 7
				11 8
				11 9
				11 10
				11 11
				11 12
				11 13
				11 14
				11 15
				11 16
				11 17
				11 18
				11 19
				11 20
				11 21
				11 22
				12 1
				12 2
				12 3
				12 4
				12 5
				12 6
				12 7
				12 8
				12 9
				12 10
				12 11
				12 12
				12 13
				12 14
				12 15
				12 16
				12 17
				12 18
				12 19
				12 20
				12 21
				13 1
				13 2
				13 3
				13 4
				13 5
				13 6
				13 7
				13 8
				13 9
				13 10
				13 11
				13 12
				13 13
				13 14
				13 15
				13 16
				13 17
				13 18
				13 19
				13 20
				14 9
				14 10
				14 11
				14 12
				14 13
				14 14
				14 15
				14 16
				14 17
				14 18
				14 19
				14 20
				14 21
				15 2
				15 3
				15 4
				15 5
				15 6
				15 7
				15 8
				15 9
				15 10
				15 11
				15 12
				15 13
				15 14
				15 15
				15 16
				15 17
				15 18
				15 19
				15 20
				15 21
				15 22
				15 23
				15 24
				16 1
				16 2
				16 3
				16 4
				16 5
				16 6
				16 7
				16 8
				16 9
				16 10
				16 11
				16 12
				16 13
				16 14
				16 15
				16 16
				16 17
				16 18
				16 19
				16 20
				16 21
				17 1
				17 2
				17 3
				17 4
				17 5
				17 6
				17 7
				17 8
				17 9
				17 10
				17 11
				17 12
				17 13
				17 14
				17 15
				17 16
				17 17
				17 18
				17 19
				17 20
				18 1
				18 2
				18 3
				18 4
				18 5
				18 6
				18 7
				18 8
				18 9
				18 10
				18 11
				18 12
				18 13
				18 14
				18 15
				18 16
				18 17
				18 18
				18 19
				18 20
				18 21
				21 1
				21 2
				21 3
				21 4
				21 5
				21 6
				21 7
				21 8
				21 9
				21 10
				21 11
				21 12
				21 13
				21 14
				21 15
				21 16
				21 17
				21 18
				21 19
				21 20
				21 21
				21 22
				21 23
				21 24
				22 1
				22 2
				22 3
				22 4
				22 5
				22 6
				22 7
				22 8
				22 9
				22 10
				22 11
				22 12
				22 13
				22 14
				22 15
				22 16
				22 17
				22 18
				22 19
				22 20
				22 21
				22 22
				22 23
				22 24
				23 1
				23 2
				23 3
				23 4
				23 5
				23 6
				23 7
				23 8
				23 9
				23 10
				23 11
				23 12
				23 13
				23 14
				23 15
				23 16
				23 17
				23 18
				23 19
				23 20
				23 21
				23 22
				23 23
				23 24
				24 1
				24 2
				24 3
				24 4
				24 5
				24 6
				24 7
				24 8
				24 9
				24 10
				24 11
				24 12
				24 13
				24 14
				24 15
				24 16
				24 17
				24 18
				24 19
				24 20
				24 21
				24 22
				24 23
				24 24
				25 1
				25 2
				25 3
				25 4
				25 5
				25 6
				25 7
				25 8
				25 9
				25 10
				25 11
				25 12
				25 13
				25 14
				25 15
				25 16
				25 17
				25 18
				25 19
				25 20
				25 21
				25 22
				25 23
				25 24
				26 1
				26 2
				26 3
				26 4
				26 5
				26 6
				26 7
				26 8
				26 9
				26 10
				26 11
				26 12
				26 13
				26 14
				26 15
				26 16
				26 17
				26 18
				26 19
				26 20
				26 21
				26 22
				26 23
				26 24
				27 1
				27 2
				27 3
				27 4
				27 5
				27 6
				27 7
				27 8
				27 9
				27 10
				27 11
				27 12
				27 13
				27 14
				27 15
				27 16
				27 17
				27 18
				27 19
				27 20
				27 21
				27 22
				27 23
				27 24
				28 1
				28 2
				28 3
				28 4
				28 5
				28 6
				28 7
				28 8
				28 9
				28 10
				28 11
				28 12
				28 13
				28 14
				28 15
				28 16
				28 17
				28 18
				28 19
				28 20
				28 21
				28 22
				28 23
				28 24
				29 1
				29 2
				29 3
				29 4
				29 5
				29 6
				29 7
				29 8
				29 9
				29 10
				29 11
				29 12
				29 13
				29 14
				29 15
				29 16
				29 17
				29 18
				29 19
				29 20
				29 21
				29 22
				29 23
				29 24
				30 1
				30 2
				30 3
				30 4
				30 5
				30 6
				30 7
				30 8
				30 9
				30 10
				30 11
				30 12
				30 13
				30 14
				30 15
				30 16
				30 17
				30 18
				30 19
				30 20
				30 21
				30 22
				30 23
				30 24
				31 1
				31 2
				31 3
				31 4
				31 5
				31 6
				31 7
				31 8
				31 9
				31 10
				31 11
				31 12
				31 13
				31 14
				31 15
				31 16
				31 17
				31 18
				31 19
				31 20
				31 21
				31 22
				31 23
				31 24
				32 1
				32 2
				32 3
				32 4
				32 5
				32 6
				32 7
				32 8
				32 9
				32 10
				32 11
				32 12
				32 13
				32 14
				32 15
				32 16
				32 17
				32 18
				32 19
				32 20
				32 21
				32 22
				32 23
				32 24
				33 1
				33 2
				33 3
				33 4
				33 5
				33 6
				33 7
				33 8
				33 9
				33 10
				33 11
				33 12
				33 13
				33 14
				33 15
				33 16
				33 17
				33 18
				33 19
				33 20
				33 21
				33 22
				33 23
				33 24
			};
			
			\addplot[
			only marks,
			mark=x,
			mark options={draw=red},
			mark size=2.5pt
			] table {
				x y
				3 1
				3 2
				3 3
				3 4
				3 5
				3 6
				3 7
				3 8
				3 9
				3 10
				3 11
				3 12
				3 13
				3 14
				3 15
				3 16
				3 17
				3 18
				3 19
				3 20
				3 21
				3 22
				3 23
				3 24
				4 1
				4 2
				4 3
				4 4
				4 5
				4 6
				4 7
				4 8
				4 9
				4 10
				4 11
				4 12
				4 13
				4 14
				4 15
				4 16
				4 17
				4 18
				4 19
				4 20
				4 21
				4 22
				4 23
				4 24
				7 1
				7 2
				7 3
				7 4
				7 5
				7 6
				7 7
				7 8
				7 9
				7 10
				7 11
				7 12
				7 13
				7 14
				7 15
				7 16
				7 17
				7 18
				7 19
				7 20
				7 21
				7 22
				7 23
				7 24
				8 1
				8 2
				8 3
				8 4
				8 5
				8 6
				8 7
				8 8
				8 9
				8 10
				8 11
				8 12
				8 13
				8 14
				8 15
				8 16
				8 17
				8 18
				8 19
				8 20
				8 21
				8 22
				8 23
				8 24
				9 1
				9 2
				9 3
				9 4
				9 5
				9 6
				9 7
				9 8
				9 9
				9 10
				9 11
				9 12
				9 13
				9 14
				9 15
				9 16
				9 17
				9 18
				9 19
				9 20
				9 21
				9 22
				9 23
				9 24
				10 1
				10 2
				10 3
				10 4
				10 5
				10 6
				10 7
				10 8
				10 9
				10 10
				10 11
				10 12
				10 13
				10 14
				10 15
				10 16
				10 17
				10 18
				10 19
				10 20
				10 21
				10 22
				10 23
				10 24
				11 1
				11 2
				11 3
				11 4
				11 5
				11 6
				11 7
				11 8
				11 9
				11 10
				11 11
				11 12
				11 13
				11 14
				11 15
				11 16
				11 17
				11 18
				11 19
				11 20
				11 21
				11 22
				11 23
				11 24
				12 1
				12 2
				12 3
				12 4
				12 5
				12 6
				12 7
				12 8
				12 9
				12 10
				12 11
				12 12
				12 13
				12 14
				12 15
				12 16
				12 17
				12 18
				12 19
				12 20
				12 21
				12 22
				12 23
				12 24
				13 1
				13 2
				13 3
				13 4
				13 5
				13 6
				13 7
				13 8
				13 9
				13 10
				13 11
				13 12
				13 13
				13 14
				13 15
				13 16
				13 17
				13 18
				13 19
				13 20
				13 21
				13 22
				13 23
				13 24
				14 11
				14 12
				14 13
				14 14
				14 15
				14 16
				14 17
				14 18
				14 19
				14 20
				14 21
				14 22
				14 23
				14 24
				15 1
				15 2
				15 3
				15 4
				15 5
				15 6
				15 7
				15 8
				15 9
				15 10
				15 11
				15 12
				15 13
				15 14
				15 15
				15 16
				15 17
				15 18
				15 19
				15 20
				15 21
				15 22
				15 23
				15 24
				16 1
				16 2
				16 5
				16 6
				16 7
				16 8
				16 9
				16 10
				16 16
				16 17
				16 18
				16 19
				16 22
				16 23
				16 24
				17 1
				17 2
				17 5
				17 6
				17 7
				17 8
				17 9
				17 10
				17 11
				17 12
				17 13
				17 16
				17 17
				17 18
				17 19
				17 22
				17 23
				17 24
				18 1
				18 2
				18 5
				18 6
				18 7
				18 8
				18 9
				18 10
				18 11
				18 12
				18 13
				18 16
				18 17
				18 18
				18 19
				18 22
				18 23
				18 24
				19 4
				19 5
				19 6
				19 7
				19 8
				19 9
				19 10
				19 11
				19 12
				19 13
				19 14
				19 15
				19 16
				19 17
				19 18
				19 19
				19 22
				19 23
				19 24
				20 4
				20 5
				20 6
				20 7
				20 8
				20 9
				20 10
				20 11
				20 12
				20 13
				20 14
				20 15
				20 16
				20 17
				20 18
				20 19
				20 22
				20 23
				20 24
				21 1
				21 2
				21 3
				21 4
				21 5
				21 6
				21 7
				21 8
				21 9
				21 10
				21 11
				21 12
				21 13
				21 14
				21 15
				21 16
				21 17
				21 18
				21 19
				21 20
				21 21
				21 22
				21 23
				21 24
				22 1
				22 2
				22 3
				22 4
				22 5
				22 6
				22 7
				22 8
				22 9
				22 10
				22 11
				22 12
				22 13
				22 14
				22 15
				22 16
				22 17
				22 18
				22 19
				22 20
				22 21
				22 22
				22 23
				22 24
				23 1
				23 2
				23 3
				23 4
				23 5
				23 6
				23 7
				23 8
				23 9
				23 10
				23 11
				23 12
				23 13
				23 14
				23 15
				23 16
				23 17
				23 18
				23 19
				23 20
				23 21
				23 22
				23 23
				23 24
				24 1
				24 2
				24 3
				24 4
				24 5
				24 6
				24 7
				24 8
				24 9
				24 10
				24 11
				24 12
				24 13
				24 14
				24 15
				24 16
				24 17
				24 18
				24 19
				24 20
				24 21
				24 22
				24 23
				24 24
				25 1
				25 2
				25 3
				25 4
				25 5
				25 6
				25 7
				25 8
				25 9
				25 10
				25 11
				25 12
				25 13
				25 14
				25 15
				25 16
				25 17
				25 18
				25 19
				25 20
				25 21
				25 22
				25 23
				25 24
				26 1
				26 2
				26 3
				26 4
				26 5
				26 6
				26 7
				26 8
				26 9
				26 10
				26 11
				26 12
				26 13
				26 14
				26 15
				26 16
				26 17
				26 18
				26 19
				26 20
				26 21
				26 22
				26 23
				26 24
				27 1
				27 2
				27 3
				27 4
				27 5
				27 6
				27 7
				27 8
				27 9
				27 10
				27 11
				27 12
				27 13
				27 14
				27 15
				27 16
				27 17
				27 18
				27 19
				27 20
				27 21
				27 22
				27 23
				27 24
				28 1
				28 2
				28 3
				28 4
				28 5
				28 6
				28 7
				28 8
				28 9
				28 10
				28 11
				28 12
				28 13
				28 14
				28 15
				28 16
				28 17
				28 18
				28 19
				28 20
				28 21
				28 22
				28 23
				28 24
				29 1
				29 2
				29 3
				29 4
				29 5
				29 6
				29 7
				29 8
				29 9
				29 10
				29 11
				29 12
				29 13
				29 14
				29 15
				29 16
				29 17
				29 18
				29 19
				29 20
				29 21
				29 22
				29 23
				29 24
				30 1
				30 2
				30 3
				30 4
				30 5
				30 6
				30 7
				30 8
				30 9
				30 10
				30 11
				30 12
				30 13
				30 14
				30 15
				30 16
				30 17
				30 18
				30 19
				30 20
				30 21
				30 22
				30 23
				30 24
				31 1
				31 2
				31 3
				31 4
				31 5
				31 6
				31 7
				31 8
				31 9
				31 10
				31 11
				31 12
				31 13
				31 14
				31 15
				31 16
				31 17
				31 18
				31 19
				31 20
				31 21
				31 22
				31 23
				31 24
				32 1
				32 2
				32 3
				32 4
				32 5
				32 6
				32 7
				32 8
				32 9
				32 10
				32 11
				32 12
				32 13
				32 14
				32 15
				32 16
				32 17
				32 18
				32 19
				32 20
				32 21
				32 22
				32 23
				32 24
				33 1
				33 2
				33 3
				33 4
				33 5
				33 6
				33 7
				33 8
				33 9
				33 10
				33 11
				33 12
				33 13
				33 14
				33 15
				33 16
				33 17
				33 18
				33 19
				33 20
				33 21
				33 22
				33 23
				33 24
			};
			
			\legend{Initial, Robust}
			
		\end{axis}
	\end{tikzpicture}
	\caption{Comparison of unit commitments for the original and robust models}
	\label{fig:commitment_compare}
\end{figure}
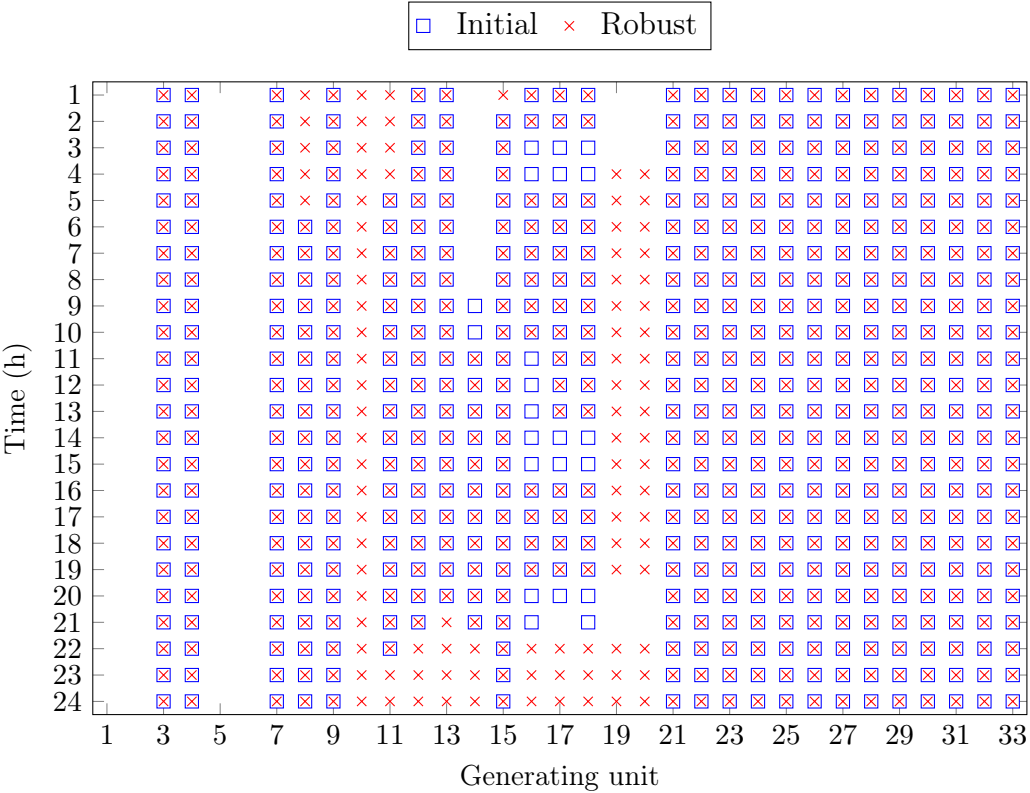

We also depict the commitment of the proposed robust model and the initial (not robust, Problem~(\ref{eq:acncuc}) without line-off) one in Figure~\ref{fig:commitment_compare}. A blue box in $(i, j)$ pertains to the initial model and indicates that the $i$-th unit is on in the $j$-th time, while the red cross mark represents the robust model.

We observe that:
\begin{enumerate}[label=\arabic*)]
	\item  The column-and-constraint generation algorithm converge in 3 iterations.
	\item  All master problems and subproblems can be solved within 160 second.
	\item  We conclude from Table~{\ref{tab:24_comp}} that the robust model provides strong protection by drastically reducing the unserved energy cost, though the robust model turns on more generating units (as Figure~\ref{fig:commitment_compare} shows) improving the commitment and production costs.
\end{enumerate}

{
\subsubsection{AC Feasibility}

We compare the proposed robust unit commitment model based on the second-order cone relaxation with two MILP-based approximations: the DC approximation and the Flat-Start AC approximation \citep{coffrin2014linear}. The detailed formulations of the DC and Flat-Start models are provided in \ref{app:dc} and~\ref{app:fs}, respectively. Table~\ref{tab:ac_feasibility} reports the computational performance of the three models within the column-and-constraint generation framework. Here, ``served cost'' denotes the generation and commitment cost, while ``unserved cost'' denotes the penalty associated with load shedding. The column ``Worst Traj.'' gives the worst-case hurricane trajectory identified by the robust model at termination, and ``Selected Traj.'' lists the trajectories (in orders) added to the master problem during the CCG iterations.

\begin{table}[!htbp]
	\centering
	\setlength{\tabcolsep}{1.1mm}
\renewcommand{\arraystretch}{1.2}
	\begin{threeparttable}
		\begin{tabular}{l|cccccc}
			\toprule[1pt]
			Model &Served Cost &Unserved Cost &Worst Traj. &Selected Traj. &\# Iter. &Time \\
			\hline
			DC         &$1.02 \times 10^{6}$ &$1.90 \times 10^{8}$ &5 &6,2,8,5 &4 &57\\
			Flat-Start &$1.02 \times 10^{6}$ &$1.90 \times 10^{8}$ &5 &6,2,8,5 &4 &106\\
			SOC        &$1.04 \times 10^{6}$ &$1.84 \times 10^{8}$ &5 &5,6,8   &3 &462\\
			\bottomrule[1pt]
		\end{tabular}
	\end{threeparttable}
	\caption{Performance comparison between DC UC, Flat-Start AC UC, and SOC AC UC on IEEE 24-bus system}
	\label{tab:ac_feasibility}
\end{table}

As shown in Table~\ref{tab:ac_feasibility}, the DC and Flat-Start approximations produce very similar robust schedules and identify the same worst-case trajectory. Both MILP-based approximations are computationally faster than the SOC relaxation, although they require one additional CCG iteration. In contrast, the SOC relaxation yields a slightly higher served cost but a lower unserved cost, suggesting that the SOC model produces a more conservative and resilient schedule under the selected worst-case trajectory.

We further evaluate the AC feasibility of the obtained solutions. Since the DC approximation does not model reactive power flow or voltage-related variables, its violation with respect to the AC nonlinear constraints is not directly reported. For the Flat-Start and SOC solutions, we evaluate the violations of two key classes of nonlinear AC constraints:
\[ c_{n,m,t}^2 + s_{n,m,t}^2 = c_{n,n,t}c_{m,m,t}, \qquad p_{n,m,t}^2 + q_{n,m,t}^2 \leq S_{n,m}^2. \]
Specifically, Table~\ref{tab:ac_feasibility_2} reports the maximum violation and the number of constraints whose violation is at least $10^{-5}$. In addition, to assess the quality of the commitment decisions under the original AC model, we fix the unit commitment decisions obtained from each approximation and locally solve the resulting continuous AC unit commitment problem using IPOPT. The corresponding AC served and unserved costs are also reported.

\begin{table}[!htbp]
	\centering
	\setlength{\tabcolsep}{1.1mm}
	\begin{threeparttable}
		\begin{tabular}{l|cccc}
			\toprule[1pt]
			Model &Max Vio &\# Vio &AC Served Cost &AC Unserved Cost \\
			\hline
			DC         &-      &-  &$1.34 \times 10^{6}$ &$1.92 \times 10^{8}$\\
			Flat-Start &0.1656 &56 &$1.17 \times 10^{6}$ &$1.93 \times 10^{8}$\\
			SOC        &0.0527 &72 &$1.04 \times 10^{6}$ &$1.84 \times 10^{8}$\\
			\bottomrule[1pt]
		\end{tabular}
	\end{threeparttable}
	\caption{AC Feasibility of DC UC, Flat-Start AC UC, and SOC AC UC on IEEE 24-bus system}
	\label{tab:ac_feasibility_2}
\end{table}

The results in Table~\ref{tab:ac_feasibility_2} indicate that the SOC relaxation provides a solution that is closer to AC feasibility in terms of the maximum nonlinear constraint violation. Although the Flat-Start approximation has fewer violated constraints, its maximum violation is substantially larger than that of the SOC relaxation. Moreover, the SOC-based commitment achieves the lowest AC unserved cost among the three models. This suggests that, despite its higher computational cost, the SOC relaxation better captures the AC network physics and leads to commitment decisions that are more reliable when evaluated under the full AC model.}

\subsection{Case study: Central Illinois 200-bus test system}

We consider 15 hurricane trajectories for the Central Illinois 200-bus test system using a similar strategy as that used for IEEE 24-bus power system. Additional details can be found in \ref{app:simulation}. 

{ We use a convergence tolerance $\epsilon:=10^{-4}$ for the column-and-constraint generation algorithm. Additionally, to solve the large-scale master problem, we adopt Alg.~\ref{Alg:oa} with $\epsilon:=10^{-4}$ (the outer-inner cutting-plane algorithm has the same convergence tolerance than the column-and-constraint generation algorithm), $\epsilon_{tol}:=10^{-5}$, $\epsilon_{par}:=\frac{1}{2}\times(10^{-5})$, and $p_{cut}:=0.55$. 

We first compare the performance between Gurobi and our outer-inner cutting-plane algorithm on solving the master problem in different iterations. The Gurobi is in default setting with a time limit as 14400 seconds. The solution times for Gurobi and for the outer-inner cutting-plane algorithm (OICPA) are provided in Table~\ref{tab:gurobi_oa_comp}. Additionally, we present both number of nonzero elements for the same master problems (M. NNZ) and their sizes (M. Size). Furthermore, we compare the objective values of best-detected solution from Gurobi and our OICPA in Gurobi Obj and OICPA Obj, respectively.}

\renewcommand*{\thefootnote}{\fnsymbol{footnote}}

\begin{table}[!htbp]
	\centering
	\setlength{\tabcolsep}{1.1mm}
\renewcommand{\arraystretch}{1.2}
	\begin{threeparttable}
		\small
		\begin{tabular}{l|cccccc}
			\toprule[1pt]
			Iteration &Gurobi time(s) &OICPA time(s) &M. NNZ &M. Size(GB) &Gurobi Obj & OICPA Obj\\
			\hline
			Iter. 0 &3792 &4113 &392075 &0.65    &$3.3177\times 10^{5}$  &$3.3168\times 10^{5}$\\
			Iter. 1 &4711 &4264 &515652 &0.72    &$3.7648\times 10^{5}$  &$3.7643\times 10^{5}$\\
			Iter. 2 &Time Limit &6587   &1011709 &1.30        &{$4.0312\times 10^{5}$}\footnote[1]{2} &$4.0259\times 10^{5}$\\
			Iter. 3 &Time Limit &11782  &1507766 &2.24        &-                      &$4.0856\times 10^{5}$\\
			\bottomrule[1pt]
		\end{tabular}
		\begin{tablenotes}
			\footnotesize
			\item[$\ast$] Gurobi stops with a relative gap about $3\%$.
		\end{tablenotes}
	\end{threeparttable}
	\caption{Performance comparison between Gurobi and Alg.~\ref{Alg:oa} on the Central Illinois 200-bus test system}
	\label{tab:gurobi_oa_comp}
\end{table}

We observe that for the original second-order-conic unit commitment (Iter. 0)) or the master problem from Iter. 1, Gurobi can solve the considered problem with a similar runtime than the outer-inner cutting-plane algorithm. However, as the problem size increases, Gurobi generally stops or breaks. Considering the problem sizes are 1.3 or 2.2 GB and our PC has 16 GB RAM, there is potentially not enough RAM for handling Branch-and-Bound tree, cut information, and approximate formulations for the mixed-integer second-order-conic program. 

{ Additionally, in Table~\ref{tab:gurobi_oa_comp}, for iterations~0 and~1, the objective values obtained by ``OICPA Obj'' are slightly better than those reported by ``Gurobi Obj''. There are two reasons for this. First, under its default settings, Gurobi does not necessarily return a solution with zero optimality gap; instead, it terminates once a feasible solution whose optimality gap is no greater than $10^{-4}$ is found. Second, when solving mixed-integer second-order-conic programs, Gurobi sometimes suffers from numerical difficulties and, occasionally, lead to suboptimal solutions. In contrast, as discussed in the last paragraph of Section~\ref{sec:oicpa}, our outer--inner cutting-plane algorithm guarantees a solution with an optimality gap below $10^{-4}$, which can therefore be better than the solution returned by Gurobi with its default settings. Moreover, the outer--inner cutting-plane algorithm decomposes the problem into mixed-integer linear programs and continuous second-order-conic programs, a process that is numerically more robust than solving mixed-integer second-order-conic programs directly. As a result, our method is more likely to produce a higher-quality feasible solution, a fact that has also been observed in~\cite{dan2}.}

We show the performance of column-and-constraint generation algorithm on the Central Illinois 200-bus test system in Tab.~\ref{tab:200bus}. The master problems is solved by Alg.~\ref{Alg:oa} and and the subproblem is solved by Mosek.

\begin{table}[!htbp]
	\centering
	\renewcommand{\arraystretch}{1.2}
	\setlength{\tabcolsep}{1.7mm}
	\begin{threeparttable}
		\begin{tabular}{l|cc|cc|c}
			\toprule[1pt]
			Iteration & LB &M. Runtime(s) &UB &S. Runtime(s) & Gap\\
			\hline
			Iter. 0 &$3.3168\times 10^{5}$ &4113  &$1.5077\times 10^{7}$ &120 &99.78\%\\
			Iter. 1 &$3.7643\times 10^{5}$ &4264  &$4.7087\times 10^{6}$ &137 &92.0\%\\
			Iter. 2 &$4.0259\times 10^{5}$ &6587  &$5.8014\times 10^{5}$ &125 &29.84\%\\
			Iter. 3 &$4.0856\times 10^{5}$ &11782 &$4.1074\times 10^{5}$ &130  &0.1\%\footnote{3}\\
			\bottomrule[1pt]
		\end{tabular}
		\begin{tablenotes}
			\footnotesize
			\item[$\dagger$] The gap of Iter. 3 is $9.85\times 10^{-5}$, which is less than the tolerance $\epsilon = 10^{-4}$.
		\end{tablenotes}
	\end{threeparttable}
	\caption{Alg.~\ref{Alg:ccg} performance on the Central Illinois 200-bus test system}
	\label{tab:200bus}
\end{table}

We observe that:
\begin{enumerate}[label=\arabic*)]
	\item  The column-and-constraint generation algorithm converge in 3 iterations.
	\item  The master problems takes significantly longer times as the number of iterations increases, while subproblems can be solved in around 130 seconds.
	\item  In Iter. 3, there is a final gap $9.85\times 10^{-5} < \epsilon$.
	\item  The comparison between commitments from the initial (not robust) model and the robust one for the scenario involving the worst hurricane is presented in Table~{\ref{tab:200_comp}}. The robust model drastically improves security by eliminating unserved energy.
\end{enumerate}

\begin{table}[!htbp]
	\centering
	\setlength{\tabcolsep}{1.7mm}
	\begin{threeparttable}
\renewcommand{\arraystretch}{1.2}
		\begin{tabular}{l|cc}
			\toprule[1pt]
			Commitment &Commitment + Served Energy Cost  &Unserved Energy Cost\\
			\hline
			Initial (Not Robust) Model &$2.9249\times 10^{5}$ &$1.2709\times 10^{7}$\\
			Robust Model & $4.1074\times 10^{5}$ &0\\
			\bottomrule[1pt]
		\end{tabular}
	\end{threeparttable}
	\caption{Initial (no robust) SCUC vs. Robust SCUC: Cost breakdown on Central Illinois 200-bus test system}
	\label{tab:200_comp}
\end{table}

\subsubsection{Solution details for outer-inner cutting-plane algorithm}

As Table~\ref{tab:gurobi_oa_comp} shows, our outer-inner cutting-plane algorithm solves 4 different master problems, and we provide some solution details for them below.

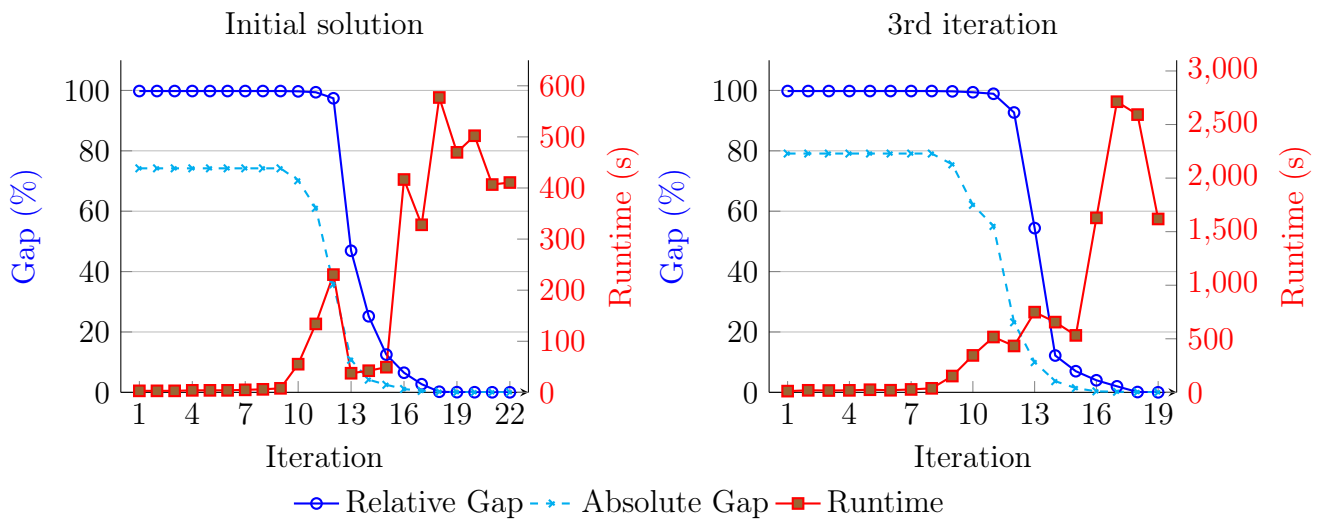
\begin{figure}[!htbp]
	\centering
	\begin{tikzpicture}
		
		\begin{axis}[
			name=plot1,
			width=7cm,
			height=6cm,
			xlabel={Iteration},
			ylabel={Gap (\%)},
			ymin=0, ymax=110,
			xtick={1,4,7,10,13,16,19,22},
			ytick={0,20,40,60,80,100},
			axis y line*=left,
			axis x line=bottom,
			ymajorgrids=true,
			enlarge x limits=0.05,
			bar width=6pt,
			tick align=inside,
			ylabel style={blue},
			every axis plot/.append style={thick},
			legend style={draw=none},
			legend to name=mylegend,
			legend columns=3,
			title={Initial solution},
			]
			\addplot+[mark=o, color=blue] coordinates {
				(1,99.8) (2,99.8) (3,99.8) (4,99.8) (5,99.8) (6,99.8) (7,99.8) (8,99.8)
				(9,99.8) (10,99.7) (11,99.4) (12,97.4) (13,46.9) (14,25.2) (15,12.5)
				(16,6.5) (17,2.7) (18,0.2) (19,0.06) (20,0.04) (21,0.02) (22,0.01)};
			\addlegendentry{Relative Gap}
			
			\addplot+[mark=x, dashed, color=cyan] coordinates {
				(1,74.2) (2,74.2) (3,74.2) (4,74.2) (5,74.2) (6,74.2) (7,74.2) (8,74.2)
				(9,74.2) (10,70.2) (11,61.1) (12,35.8) (13,10.77) (14,4.2) (15,2.5)
				(16,1.1) (17,0.48) (18,0.21) (19,0.06) (20,0.04) (21,0.02) (22,0.01)};
			\addlegendentry{Absolute Gap}
			
			\addplot+[mark=square*, color=red] coordinates {
				(1,3/6) (2,3/6) (3,3/6) (4,4/6) (5,4/6) (6,4/6) (7,5/6) (8,6/6) (9,8/6)
				(10,56/6) (11,136/6) (12,234/6) (13,38/6) (14,43/6) (15,50/6) (16,423/6)
				(17,333/6) (18,586/6) (19,477/6) (20,510/6) (21,413/6) (22,417/6)};
			\addlegendentry{Runtime}
		\end{axis}
		
		\begin{axis}[
			width=7cm,
			height=6cm,
			at={(plot1.south west)},
			anchor=south west,
			axis y line*=right,
			axis x line=none,
			ymin=0, ymax=650,
			ylabel={Runtime (s)},
			ylabel style={color=red},
			yticklabel style={color=red},
			ytick={0,100,200,300,400,500,600},
			]
			\addplot[color=white] coordinates {(0,0)};
		\end{axis}
		
		\begin{axis}[
			name=plot2,
			at={(plot1.east)},
			anchor=west,         
			xshift=3.2cm, 
			width=7cm,
			height=6cm,
			xlabel={Iteration},
			ylabel={Gap (\%)},
			ymin=0, ymax=110,
			xtick={1,4,7,10,13,16,19,22},
			ytick={0,20,40,60,80,100},
			axis y line*=left,
			axis x line=bottom,
			ymajorgrids=true,
			enlarge x limits=0.05,
			bar width=6pt,
			tick align=inside,
			ylabel style={blue},
			every axis plot/.append style={thick},
			title={3rd iteration},
			]
			\addplot+[mark=o, color=blue] coordinates {
				(1,99.8) (2,99.8) (3,99.8) (4,99.8) (5,99.8) (6,99.8) (7,99.8) (8,99.8)
				(9,99.7) (10,99.4) (11,98.9) (12,92.7) (13,54.4) (14,12.2) (15,7)
				(16,4) (17,2) (18,0.1) (19,0.01)};
			\addplot+[mark=x, dashed, color=cyan] coordinates {
				(1,79.1) (2,79.1) (3,79.1) (4,79.1) (5,79.1) (6,79.1) (7,79.1) (8,79.1)
				(9,75.5) (10,62.1) (11,55.1) (12,23.2) (13,10) (14,3.7) (15,1.4)
				(16,0.33) (17,0.2) (18,0.1) (19,0.01)};
			\addplot+[mark=square*, color=red] coordinates {
				(1,12.26/30) (2,20.26/30) (3,18.4/30) (4,20/30) (5,26/30) (6,21/30)
				(7,29/30) (8,40/30) (9,162/30) (10,367/30) (11,550/30) (12,461/30)
				(13,797/30) (14,698/30) (15,566/30) (16,1734/30) (17,2888/30)
				(18,2760/30) (19,1723/30)};
		\end{axis}
		
		\begin{axis}[
			width=7cm,
			height=6cm,
			at={(plot2.south west)},
			anchor=south west,
			axis y line*=right,
			axis x line=none,
			ymin=0, ymax=3100,
			ylabel={Runtime (s)},
			ylabel style={color=red},
			yticklabel style={color=red},
			ytick={0,500,1000,1500,2000,2500,3000},
			]
			\addplot[color=white] coordinates {(0,0)};
		\end{axis}
		
	\end{tikzpicture}
	
	\ref{mylegend}
	
	\caption{Outer-inner cutting-plane algorithm performance: (Left) initial solution, (Right) 3rd iteration}
	\label{fig:two_sub}
\end{figure}

In Fig.~\ref{fig:two_sub}, we present the following characteristics for the initial and final iterations:
\begin{enumerate}[label=\arabic*)]
	\item The blue line with mark ``o'' is the relative gap, i.e., $\frac{\mathrm{UB} - \mathrm{LB}}{\mathrm{UB}}$. The relative gap illustrate the converge evolution.
	\item The cyan line with mark ``x'' is the absolute gap, which is $\frac{obj^* - \mathrm{LB}}{obj^*}$, where $obj^*$ is the objective value of the final solution. The absolute gap provides the distance between the current solution to the final solution.
	\item The red line with the square mark is the Gurobi runtime in each iteration. We adopt a warm start from the last solution, and sometimes we get a runtime reduction, e.g. iteration 13, 14, 15 in Fig.~\ref{fig:two_sub} (left).
\end{enumerate} 

The left y-axis, e.g. Gap (\%), is for relative and absolute gaps, while the right one, e.g., Runtime (s), is for runtime.

{
\subsubsection{AC Feasibility}

We compare the proposed SOC-based robust model with the DC and Flat-Start approximations on the Central Illinois 200-bus test system. The results are reported in Tables~\ref{tab:ac_feasibility_3} and~\ref{tab:ac_feasibility_4}.

\begin{table}[!htbp]
	\centering
	\setlength{\tabcolsep}{1.1mm}
	\renewcommand{\arraystretch}{1.2}
	\begin{threeparttable}
		\begin{tabular}{l|cccccc}
			\toprule[1pt]
			Model &Served Cost &Unserved Cost &Worst Traj. &Selected Traj. &\# Iter. &Time \\
			\hline
			DC         &$3.95\times 10^{5}$ &$9.80\times 10^{5}$ &1 &3,6,1,14,2 &5 &2263\\
			Flat-Start &$3.95\times 10^{5}$ &$9.80\times 10^{5}$ &1 &1,2,3,14,6 &5 &3856\\
			SOC        &$4.11\times 10^{5}$ &0                   &2 &14, 6, 2   &3 &27258\\
			\bottomrule[1pt]
		\end{tabular}
	\end{threeparttable}
	\caption{Performance comparison between DC UC, Flat-Start AC UC, and SOC AC UC on the Central Illinois 200-bus test system}
	\label{tab:ac_feasibility_3}
\end{table}

Table~\ref{tab:ac_feasibility_3} shows that the DC and Flat-Start approximations produce nearly identical objective components and identify the same worst-case trajectory. Both models are faster than the SOC relaxation, as expected, since they lead to MILP approximations. However, both approximations result in a positive unserved cost of $9.84 \times 10^5$. In contrast, the SOC relaxation obtains a schedule with zero unserved cost, although at the expense of a higher served cost and a longer solution time. This indicates that the SOC relaxation produces a more conservative commitment and dispatch decision, but one that is substantially more robust against the considered worst-case events.

\begin{table}[!htbp]
	\centering
	\setlength{\tabcolsep}{1.1mm}
	\begin{threeparttable}
\renewcommand{\arraystretch}{1.2}
		\begin{tabular}{l|cccc}
			\toprule[1pt]
			Model &Max Vio &\# Vio &AC Served Cost &AC Unserved Cost \\
			\hline
			DC &- &- &$3.99\times10^{5}$ &$9.84\times10^{5}$\\
			Flat-Start &6.6005 &8502 &$3.99\times10^{5}$ &$9.84\times10^{5}$\\
			SOC &$4.26\times10^{-5}$ &372 &$4.15\times10^{5}$ &0\\
			\bottomrule[1pt]
		\end{tabular}
	\end{threeparttable}
	\caption{AC Feasibility of DC UC, Flat-Start AC UC, and SOC AC UC on the Central Illinois 200-bus test system}
	\label{tab:ac_feasibility_4}
\end{table}

The AC feasibility assessment in Table~\ref{tab:ac_feasibility_4} further highlights the difference between the MILP approximations and the SOC relaxation. The Flat-Start approximation has a large maximum violation, $6.6005$, and a large number of violated nonlinear AC constraints. Although the SOC solution still has some small violations above the reporting tolerance, these violations are numerically minor compared with those of the Flat-Start solution. Additionally, the DC and Flat-Start commitments both lead to an AC unserved cost of $9.84 \times 10^5$, while the SOC-based commitment maintains zero unserved cost under the full AC model. These results suggest that the SOC relaxation better preserves the AC network physics and yields commitment decisions that are more compatible with the full AC formulation.}

\section{Conclusions}\label{sec:conclusion}

{
	This paper proposes a robust optimization framework for scheduling electricity generation units under extreme weather events. The problem is formulated as a tri-level min--max--min model, explicitly capturing the interaction between operational decisions and worst-case system disruptions.
	\par
	To ensure tractability, the framework integrates second-order conic relaxations with duality-based reformulations, enabling an efficient reduction of the original problem structure. Building on this reformulation, we employ a column-and-constraint generation algorithm to solve the resulting model. In particular, the mixed-integer second-order conic master problem is handled through a tailored outer-approximation procedure, which significantly improves computational performance.
	\par
	The proposed methodology is first illustrated using a simple example to clarify its main components, and its effectiveness is subsequently demonstrated through a case study that highlights its practical applicability.
	\par
	Future research may focus on developing tighter power flow representations, enhancing decomposition and acceleration techniques, and exploring alternative uncertainty models that more accurately capture the spatial and temporal dynamics of extreme weather events.
}

\appendix
\section{Details for hurricane trajectory simulations}\label{app:simulation}

We have following 8 hurricane trajectories for the 24-bus system. The lines reported below are the ones disabled.
\begin{enumerate}
	\item (3, 24), (9, 11), (10, 12)
	\item (3, 1), (2, 4), (2, 6), (7, 8)
	\item (2, 6), (3, 1), (4, 9), (5, 10), (8, 9)
	\item (2, 6), (3, 9), (3, 24), (4, 9), (5, 10), (8, 9), (8, 10)
	\item (11, 13), (11, 14), (12, 13), (12, 23), (15, 24)
	\item (12, 23), (13, 23), (14, 16), (15, 16), (15, 21)
	\item (16, 17), (16, 19), (17, 22)
	\item (15, 21), (16, 17), (17, 22), (21, 22)
\end{enumerate}

We have following 15 hurricane trajectories for the 200-bus system. The lines reported below are the ones disabled.
\begin{enumerate}
	\item (25, 199), (171, 195)
	\item (25, 199), (81, 82), (59, 119)
	\item (93, 191), (100, 184), (14, 121), (97, 186), (109, 186)
	\item (93, 191), (63, 184), (57, 159)
	\item (34, 137), (14, 149), (97, 186), (109, 186)
	\item (85, 120), (42, 44), (58, 95), (45, 187), (46, 122), (81, 178), (25, 64)
	\item (60, 134), (128, 133), (43, 132), (144, 162), (59, 119)
	\item (93, 191), (100, 184), (141, 121), (42, 44), (85, 120)
	\item (25, 199), (81, 178), (46, 122), (45, 187), (58, 177), (83, 146), (60, 186)
	\item (34, 54), (14, 149), (58, 95), (31, 192), (60, 134)
	\item (25, 199), (81, 178), (46, 122), (45, 181), (31, 192), (60, 134)
	\item (34, 54), (14, 149), (83, 186), (60, 186)
	\item (93, 191), (63, 184), (160, 181), (39, 85)
	\item (93, 191), (100, 184), (14, 121), (44, 200), (17, 109)
	\item (34, 54), (14, 15), (14, 121), (58, 95), (31, 192), (60, 186)
\end{enumerate}

{
\section{Formulation for DC Unit Commitment}\label{app:dc}

\begin{subequations}
	\begin{align}
		\min\ &\sum_{t\in\mathcal{T}}\sum_{g\in\mathcal{G}} C_{g}^\mathrm{F}u_{g,t}+C_{g}^\mathrm{SU}y_{g,t}+C_{g}^\mathrm{SD}z_{g,t}
		+\sum_{t\in\mathcal{T}}\left[\sum_{g\in\mathcal{G}}C_g^\mathrm{V} p_{g,t}+\sum_{n\in\mathcal{N}} C^\mathrm{U}p_{n,t}^\mathrm{U} \right]\\
		\mbox{s.t.}\ &(\ref{eq:p_logic1}) - (\ref{eq:p_bal1})\notag\\
		&p_{n,m,t} = -B_{n,m}(\theta_{n,t} - \theta_{m,t}),\quad \forall t\in\mathcal{T}, \forall (n,m)\in \mathcal{E}\label{eq:dc_p}\\
		&-\frac{\pi}{6} \leq \theta_{n,t} - \theta_{m,t} \leq \frac{\pi}{6},\quad \forall t\in\mathcal{T}, \forall (n,m)\in \mathcal{E}\\
		&-S_{n,m} \leq p_{n,m,t} \leq S_{n,m}, \quad \forall t\in\mathcal{T}, \forall (n,m)\in \mathcal{E} \label{eq:dc_cap}
	\end{align}
\end{subequations}

\section{Formulation for Flat-Start AC Unit Commitment}\label{app:fs}

\begin{subequations}
	\begin{align}
		\min\ &\sum_{t\in\mathcal{T}}\sum_{g\in\mathcal{G}} C_{g}^\mathrm{F}u_{g,t}+C_{g}^\mathrm{SU}y_{g,t}+C_{g}^\mathrm{SD}z_{g,t}
		+\sum_{t\in\mathcal{T}}\left[\sum_{g\in\mathcal{G}}C_g^\mathrm{V} p_{g,t}+\sum_{n\in\mathcal{N}} C^\mathrm{U}(p_{n,t}^\mathrm{U} + q_{n,t}^\mathrm{U})\right]\\
		\mbox{s.t.}\ &(\ref{eq:p_logic1}) - (\ref{eq:p_bal2})\notag\\
		&p_{n,m,t} = -G_{n,m}(v_{n,t} - v_{m,t}) - B_{n,m}(\theta_{n,t} - \theta_{m,t}),\quad \forall t\in\mathcal{T}, \forall (n,m)\in \mathcal{E}\label{eq:flat_p}\\
		&q_{n,m,t} = -B_{n,m}(v_{n,t} - v_{m,t}) + G_{n,m}(\theta_{n,t} - \theta_{m,t}),\quad \forall t\in\mathcal{T}, \forall (n,m)\in \mathcal{E}\label{eq:flat_q}\\
		&V_{n}^{\min} - 1 \leq v_{n,t} \leq V_{n}^{\max} - 1,\quad \forall t\in \mathcal{T}, n\in \mathcal{N}\\
		&-\frac{\pi}{6} \leq \theta_{n,t} - \theta_{m,t} \leq \frac{\pi}{6},\quad \forall t\in\mathcal{T}, \forall (n,m)\in \mathcal{E}\\
		&-S_{n,m} \leq p_{n,m,t}, q_{n,m,t} \leq S_{n,m}, \quad \forall t\in\mathcal{T}, \forall (n,m)\in \mathcal{E} \label{eq:cap_1}\\
		&-\sqrt{2}S_{n,m} \leq p_{n,m,t} + q_{n,m,t} \leq \sqrt{2}S_{n,m}, \quad \forall t\in\mathcal{T}, \forall (n,m)\in \mathcal{E} \label{eq:cap_2}
	\end{align}
\end{subequations}
We introduce $v_{n,t} = 1 - |V_{n,t}|$, where $V_{n,t}$ is the complex voltage of node $n$ in $t$, and $\theta_{n, t}$ is the angle of node $n$ in $t$. Based on the flat-start approximation proposed by \citep{coffrin2014linear}, the active and reactive power flows between node $n$ and $m$ in time $t$ are represented by (\ref{eq:flat_p}) and (\ref{eq:flat_q}) respectively. We linearized the line capacity constraint $p_{n,m,t}^2 + q_{n,m,t}^2 \leq S_{n,m}^2$ by (\ref{eq:cap_1}) and (\ref{eq:cap_2}), where (\ref{eq:cap_2}) is from the Cauchy-Schwarz inequality $p_{n,m,t} + q_{n,m,t} \leq \sqrt{2}\sqrt{p_{n,m,t}^2 + q_{n,m,t}^2} \leq \sqrt{2} S_{n,m}$.

When test the AC feasibility, we calculate $V_{n,t} = 1 + v_{n,t}$,  $c_{n,n,t} = |V_{n,t}|^2$, $c_{n,m,t} = |V_{n,t}| |V_{m,t}| \cos(\theta_{n,t} - \theta_{m,t})$, and $s_{n,m,t} = -|V_{n,t}| |V_{m,t}| \sin(\theta_{n,t} - \theta_{m,t})$.

}

\section*{Acknowledgment}

The work reported in this paper has been partially supported by the Advanced Grid Modeling Program of the Office of Electricity, U.S. Department of Energy, through Argonne National Laboratory.

\bibliographystyle{elsarticle-harv}  
\bibliography{ref}

\end{document}